\documentclass{amsart}

\usepackage[all]{xy}   

\theoremstyle{plain}
\newtheorem{lem}{Lemma}[section]
\newtheorem{cor}[lem]{Corollary}
\newtheorem{prop}[lem]{Proposition}
\newtheorem{thm}[lem]{Theorem}

\theoremstyle{definition}
\newtheorem{ex}[lem]{Example}
\newtheorem{exs}[lem]{Examples}
\newtheorem{rem}[lem]{Remark}
\newtheorem{dfn}[lem]{Definition}

\newcommand{\Z}{\mathbb{Z}}               
\newcommand{\Q}{\mathbb{Q}}              
\newcommand{\al}{\alpha}       
\newcommand{\be}{\beta}        
\newcommand{\ga}{\gamma}   
\newcommand{\de}{\delta}       
\newcommand{\ep}{\epsilon}   
\newcommand{\la}{\lambda}    
\newcommand{\La}{\Lambda}  
\newcommand{\De}{\Delta}       
\newcommand{\id}{\mathrm{id}}      
\newcommand{\Gg}{\mathcal{G}}           
\newcommand{\SA}{\mathcal{Z}}     
\newcommand{\SC}{\mathfrak{S}}   

\newcommand{\Fp}{\mathbb{F}_p}   
\newcommand{\CH}{\mathop{\rm CH}}   

\newcommand{\hh}{\mathtt{h}}     

\newcommand{\SBM}{S\text{-}\mathsf{BMod}_R}  
\newcommand{\RBM}{R\text{-}\mathsf{BMod}_R}  

\newcommand{\hotimes}{\mathop{\widehat{\otimes}}}  
\newcommand{\botimes}{\mathop{\otimes}}   

\newcommand{\hunit}{\hat{\mathfrak{1}}}   
\newcommand{\bunit}{\mathfrak{1}}    

\def\compactstyle{%
\def\objectstyle{%
    \let\oldotimes\otimes%
    \def\otimes{\mathop{\raisebox{0.4ex}[1.25ex][0pc]{%
    \({\scriptstyle\oldotimes}\)%
    }}\limits}
    \def\SWS{S{\rule[-1ex]{0pc}{0pc}\raisebox{0.4ex}[0pc][0pc]{\hspace{0.15em}%
    \(\mathop{\scriptstyle\oldotimes}\limits_{\!\!\!\!\!\!
    S^{\hspace{-0.1em}\scriptscriptstyle W}\!\!\!\!\!\!}\)%
    \hspace{0.15em}}} S}%
    \def\SRS{S{\rule[-1ex]{0pc}{0pc}\raisebox{0.4ex}[0pc][0pc]{\hspace{0.15em}%
    \(\mathop{\scriptstyle\oldotimes}\limits_{\!\!\!\!\!\!
    R\!\!\!\!\!\!}\)%
    \hspace{0.15em}}} S}%
    \def\hotimes{{\raisebox{0.25ex}[1.2ex][0pc]{\(\mathop{\scriptstyle\widehat{\oldotimes}}\)}}}%
    \def\botimes{{\raisebox{0.25ex}[1.2ex][0pc]{\(\mathop{\scriptstyle{\oldotimes}}\)}}}%
    \let\oldwedge\wedge
    \def\wedge{{\oldwedge}}
    \let\oldvee\vee
    \def\vee{{\oldvee}}}
\def\labelstyle{\scriptstyle%
    \let\oldotimes\otimes%
    \def\otimes{{\raisebox{0.1ex}{%
    \({\scriptscriptstyle\oldotimes}\)%
    }}}
    \def\hotimes{{\raisebox{0.1ex}{%
    \({\scriptscriptstyle\widehat{\oldotimes}}\)%
    }}}
}}
\begin{document}

\title[Structure algebras, Hopf algebroids and oriented cohomology]{Structure algebras, Hopf algebroids and oriented cohomology of a group}

\author[M.~Lanini]{Martina Lanini}
\address[Martina Lanini]{Dipartimento di Matematica, University of Rome Tor Vergata, Via della Ricerca Scientifica 1, 00133, Rome, Italy}
\email{lanini@mat.uniroma2.it}
\urladdr{https://sites.google.com/site/martinalanini5/home}

\author[R.~Xiong]{Rui Xiong}
\address[Rui Xiong]{Department of Mathematics and Statistics, University of Ottawa, 150 Louis-Pasteur, Ottawa, ON, K1N 6N5, Canada}
\email{rxiong@uottawa.ca}
\urladdr{https://arxiv.org/a/0000-0001-8025-0795.html}

\author[K.~Zainoulline]{Kirill Zainoulline}
\address[Kirill Zainoulline]{Department of Mathematics and Statistics, University of Ottawa, 150 Louis-Pasteur, Ottawa, ON, K1N 6N5, Canada}
\email{kirill@uottawa.ca}
\urladdr{http://mysite.science.uottawa.ca/kzaynull/}

\subjclass[2010]{14F43, 14M15, 14L30,  18M50, 16T05}
\keywords{equivariant oriented cohomology, structure algebra, cohomology of a group, duoidal category}

\begin{abstract}  We prove that the structure algebra of a Bruhat moment graph of a finite real root system is a Hopf algebroid with respect to the Hecke and the Weyl actions. 
We introduce new techniques (reconstruction and push-forward formula of a product, twisted coproduct, double quotients of bimodules) and apply them together with our main result to linear algebraic groups, to generalized Schubert calculus, to combinatorics of Coxeter groups and finite real root systems.
As for groups, it implies that the natural Hopf-algebra structure on the algebraic oriented cohomology $\hh(G)$ of Levine-Morel of a split semi-simple linear algebraic group $G$ can be lifted to a `bi-Hopf' structure on the $T$-equivariant algebraic oriented cohomology of the complete flag variety.  As for Schubert calculus, we prove several new identities involving (double) generalized equivariant Schubert classes. As for finite real root systems, we compute the Hopf-algebra structure of `virtual cohomology' of dihedral groups $I_2(p)$, where $p$ is an odd prime. 
\end{abstract}

\maketitle

\tableofcontents


\section{Introduction}
The Chow ring $\CH(G)$ (or a singular cohomology) with integer coefficients of a split semi-simple linear algebraic group $G$ is the most celebrated geometric invariant in the theory of linear algebraic groups. Starting from pioneering works by Grothendieck and Borel, it has been studied for decades and computed for all simple groups (for overview see e.g. Kac \cite{Ka85}, Duan-Zhao \cite{DZ14}). The major difficulty here is to understand the torsion part of $\CH(G)$. Classical arguments involving the Grothendieck torsion index \cite{Gr58} show that for all simply-connected groups of type $A$ and $C$ there is no torsion, i.e., $\CH(SL_n)\simeq \CH(Sp_{2n})\simeq \Z$. It is comparably easy to compute it if $G$ has only one torsion (bad) prime $p$. For instance, $\CH(PGL_p)\simeq \Z[x]/(px,x^p)$, where $x$ is a generator of degree $1$. However, in the mixed torsion prime case the complexity of the answer rises enormously (we refer to \cite{DZ14} for examples of computations for exceptional groups). 

If we consider another cohomological-type functor $\hh$ -- the so called oriented (generalized) algebraic cohomology theory in the sense of Levine-Morel \cite{LM07},
e.g., connective $K$-theory,  hyperbolic cohomology of \cite{LZZ}, various Morava $K$-theories or the algebraic cobordism $\Omega$ of Levine-Morel (the universal such theory) -- then the answer for $\hh(G)$ remains widely open. (The only case where the situation is well-understood is the $K$-theory.) Already for rank 2 groups, only partial computations of algebraic cobordism are known (see e.g. \cite{Ya05}). Observe that the coefficient ring of $\Omega$ is not even $\Z$ but the Lazard ring (a polynomial ring in infinitely many variables).

One important feature of such $\hh(G)$ is the fact that it has a natural coproduct structure arising from the geometric multiplication $G\times G\to G$ and the inverse $G\to G$ which turns $\hh(G)$ into a Hopf algebra. This substantially reduces the computations as one can use the Hopf-Borel theorem or other structural theorems about Hopf algebras.
For instance, the Chow ring at a torsion prime splits as 
\[\CH(G;\Fp)\simeq \bigotimes_i \Fp[x_i]/(x_i^{p^{k_i}}),\] 
which reduces the problem to determining the exponents $k_i$ given by $p$-exceptional degrees of Kac~\cite{Ka85}. 
Observe that the numerical invariants $k_i$s play an important role in the study of motivic decomposition types of various versal flag varieties \cite{PSZ}, in the geometric theory of quadratic forms in the theory of canonical dimension of linear algebraic groups \cite{Za08} and other areas. 

The standard approach in computing $\CH(G)$  (and $\hh(G)$ in general) goes back to Grothendieck and relies on the use of the topological fibration 
$G\to G/B\to \mathcal{B} B$ (here $B$ is a Borel subgroup of $G$, $\mathcal{B}B$ is its classifying space, and $G/B$ is the variety of complete flags) and the induced isomorphism \[\CH(G)\simeq \CH(G/B)/(im_{>0} c),\] where $c\colon \CH(\mathcal{B} B)\to \CH(G/B)$ is the characteristic map and $(im_{>0} c)$ is the characteristic ideal (the generalization of this fact for an arbitrary $\hh$ can be found in \cite{GZ12}). 
However, via this approach it is pretty hard to identify the generators and relations (e.g. in terms of Schubert classes) and to exploit the Hopf-structure, as $\hh(G/B)$ is not even a Hopf algebra.  

In the present paper we overcome this difficulty by generalizing this approach 
to the context of structure algebras of finite real root systems and formal group laws. 
Our main object of study is an $S$-bimodule $\SA$ defined in purely algebraic terms from a root system $\Phi$, an isogeny lattice $\La$, and a formal group ring $S$. 
The bimodule $\SA$ was first introduced and studied in \cite{DLZ} 
as a (generalized) oriented analogue of the ring of global sections of the structure sheaf on a Bruhat moment graph.
If $\Phi$ is crystallographic,  $\La$ is an isogeny lattice, $G$ is the associated group and $T$ is its split maximal torus, then \[S\simeq \hh_T(pt)\quad\text{ and }\quad \SA\simeq \hh_T(G/B)\] are the equivariant coefficient ring and oriented cohomology respectively  \cite{KK86,KK90,CZZ}.
However, for a non-crystallographic $\Phi$ the ring $\SA$ while still well-defined does not have any geometric counterpart, 
and, hence, may be thought of as a virtual equivariant cohomology.

Our main result (Theorem~\ref{thm:bimonoid}) says that the structure algebra $\SA$ 
possesses some traces of the Hopf structure that is the so called Hopf algebroid structure
(or `bi-Hopf' structure) in the language of 2-monoidal (duoidal) categories. 
Roughly speaking, it is given by Hopf-similar diagrammatic relations involving two different tensor products $\otimes$ and $\hat\otimes$ instead of one on the duoidal category of $S$-bimodules: the first tensor product is given by  the tautological tensor product $s_1xs_2\otimes y=x\otimes s_1ys_2$ and the second one is the standard bimodule tensor product $xs\hat\otimes y=x\hat\otimes sy$. 

Since after taking the double quotient $\hat\SA=\SA/(I\SA+\SA I)$ (here $I$ is the augmentation ideal of $S$) both $\otimes$ and $\hat\otimes$ coincide,
we immediately obtain that $\hat\SA$ has a natural structure of a Hopf algebra (Corollary~\ref{cor:Hopfalg}). In other words, the `bi-Hopf' structure reduces to the usual Hopf structure. Observe that in the crystallographic case we have $\hat\SA=\hh(G)$, and for non-crystallographic $\Phi$ 
the Hopf algebra $\hat\SA$ produces an interesting invariant of $\Phi$ that can be viewed as a virtual cohomology `$\CH(G(\Phi))$' of a virtual group `$G(\Phi)$' of $\Phi$. 

It would be interesting to see whether our results can be extended to Kac-Moody root systems and Kac-Moody groups,
or to other examples of categories of bimodules appearing in geometric representation theory.

\vspace{3mm}

\paragraph{\it Organization of the paper}
To stress that in our proofs and constructions we only rely on purely algebraic, combinatorial and categorical arguments, 
we intentionally put all geometric facts, references and motivation (if needed) in the beginning of each section, keeping the rest `free of any geometry'.

We start with preliminaries in Section~\ref{sec:momgr}. These include root systems, root and weight lattices, formal group laws,
formal group rings, push-pull elements, the structure algebra $\SA$, Hecke and the Weyl-action on it. We finish this section by introducing the characteristic and Borel maps
and discuss its basic properties. In the next section~\ref{sec:intmult} we define `twisted Schubert classes' and the `push-forward to a point' element. We then introduce the dual basis and discuss various intersection multiplicities. In Section~\ref{sec:pushforw} we prove our first main result -- the push-forward formula for a product of characteristic classes (Theorem~\ref{thm:pushfres}). As a corollary, we obtain the reconstruction formula for double Schubert elements (Corollary~\ref{cor:recform}).
Next, in Section ~\ref{sec:coprod}, we construct the coproduct for the structure algebra $\SA$ with respect to the tensor product of bimodules $\hat\otimes$.
This is the key technical (and the most challenging) step used in the proof of our main result. In Section~\ref{sec:duoidal} we prove our main result (Theorem~\ref{thm:bimonoid}). We also show that the double quotient of the structure algebra has a natural structure of a bialgebra (Corollary~\ref{cor:bialgeb}) and a Hopf algebra (Corollary~\ref{cor:Hopfalg}). In the last section, we discuss an example of the dihedral non-crystallographic root system
$I_2(p)$, where $p$ is an odd prime and show how to compute its virtual cohomology $\hat\SA$ and the coproduct.

\vspace{3mm}

\paragraph{\it Acknowledgements} R.X. and K.Z. were partially supported by NSERC Discovery grant RGPIN-2022-03060. M.L acknowledges the PRIN2017 CUP E8419000480006, the Fondi di Ricerca Scientifica d'Ateneo 2021 CUP E83C22001680005, and the MUR Excellence Department Project MatMod@TOV awarded to the Department of Mathematics, University of Rome Tor Vergata.


\section{Preliminaries}\label{sec:momgr}


In this paper as an input we use two ingredients of purely algebraic nature: a finite real root system $\Phi$ together with an `isogeny' root lattice $\La$, and a 1-dimensional commutative formal group law $F$ over its coefficient ring $R$. The first ingredient (in the crystallographic case) gives rise to the semi-simple split linear algebraic group $G$, associated variety of complete flags $G/B$,  their usual and $T$-equivariant cohomology (where $T$ is a split maximal torus of $G$), and more generally (also in the non-crystallographic case) to the theory of structure algebras $\SA$ of moment graphs $\Gg$. The second ingredient allows us to extend naturally (or rather to deform) all these constructions to the context of an arbitrary (generalized) oriented cohomology theory $\hh(\text{-})$ and its $T$-equivariant analogue $\hh_T(\text{-})$. 
We briefly recall both inputs in~\ref{int:roots} and~\ref{int:fgrs} also introducing and fixing the needed notation. 

Merging root lattices and formal group laws leads us to the notion of a formal group ring $S$ (see~\ref{int:fgring}). It has been first introduced in~\cite[\S2]{CPZ}, and it plays a central role in computations of (generalized) oriented cohomology of projective homogeneous varieties. Indeed, as was shown in \cite[\S3]{CZZ} it coincides with the equivariant coefficient ring $\hh_T(pt)$ (completed equivariant cohomology of a point). 

In the next step (see~\ref{int:pushpull}), starting from the formal group ring $S$ we introduce several crucial constructions (of a twisted group algebra $S_W$ and its localization $Q_W$, of push-pull elements $Y_i$, and of the Demazure algebra $\mathbb{D}$ which are main technical tools in the localization theory. A key object of this paper -- the notion of a structure algebra $\SA$ is introduced in~\ref{int:struct}. Here we follow the extended definition of \cite[\S6]{DLZ}. Observe that in the crystallographic (geometric) case $\SA$ can be identified with the equivariant oriented cohomology $\hh_T(G/B)$.   We finalize this section by recalling basic properties of the Hecke '$\bullet$-' and of the Weyl '$\odot$-' actions on structure algebras, of the respective characteristic $c$ and of the Borel map $\rho$. Here we follow mostly \cite[\S3]{LZZ}.


\subsection{Root systems, root and weight lattices}\label{int:roots}
Let $\Phi$ be a finite real root system in the Euclidean space $\mathbb{R}^N$ together with a subset of simple roots $\Pi$ and the decomposition $\Phi=\Phi_+\amalg \Phi_-$ into positive and negative roots as in \cite[I,\S1]{Hu90}. We denote by $\mathcal{O}$ its coefficient ring (the smallest subring of $\mathbb{R}$ containing all coefficients expressing roots as linear combinations of simple roots). Observe that $\mathcal{O}=\Z$ in  crystallographic cases, and $\mathcal{O}=\Z(2\cos(\tfrac{\pi}{m}))$, $m=5$ or $m\ge 7$ in non-crystallographic cases.

We denote by $\La_r$ the corresponding {\it root lattice} of $\Phi$ (the $\mathcal{O}$-linear span of $\Pi$). Let $K$ be the field of fractions of $\mathcal{O}$ (it is a subfield of $\mathbb{R}$). We define the {\it weight lattice} of $\Phi$ to be the $\mathcal{O}$-module
\[
\La_w=\{\la\in \La_r\otimes_\mathcal{O} K\mid \al^\vee(\la)\in \mathcal{O} \text{ for all }\al\in \Phi_+\},
\]
where $\al^\vee$  denotes the respective coroot in $(\mathbb{R}^{N})^*$, that is $\al^\vee(\la)=\tfrac{2(\al,\la)}{(\al,\al)}$.
If the root system $\Phi$ is normalized, i.e., all simple roots have the same length 1 (e.g. in all non-crystallographic cases), then we have 
\[
\al^\vee(\la)\in \mathcal{O}\text{ for all }\al\in \Phi_+\quad \Longleftrightarrow\quad \al_j^\vee(\la)\in \mathcal{O}\text{ for all }\al_j\in\Pi.
\]
So in the normalized case $\La_w=M^{-1}(\La_r)$ is a finitely generated free $\mathcal{O}$-module, where $M=(\al_j^\vee(\al_i))_{ij}$ is the Cartan (Schl\"afli) matrix. Observe that by definition $\La_r\otimes_\mathcal{O} K=\La_w\otimes_\mathcal{O} K$, and $\La_w/\La_r$ is a finite abelian group.

\begin{ex}\label{ex:dihedral} 
Consider a {\em dihedral} root system $\Phi$ of type $I_2(p)$, $p$ is an odd prime, corresponding to the symmetry group of a regular $2p$-gon. (Observe that for $p=3$ it gives a crystallographic root system of type $A_2$.) On the complex plain $\Phi=\{1,\xi,\ldots,\xi^{2p-1}\}$, where $\xi$ is a primitive $2p$-th root of unity, with positive roots $\Phi_+=\{1,\ldots,\xi^{p-1}\}$ and simple roots $\Pi=\{\al_1=1, \al_2=\xi^{p-1}\}$. Its coefficient ring $\mathcal{O}$ is the ring of algebraic integers $\Z(\ga)$ in the formally real cyclotomic field $\Q(\ga)$, where $\ga=\xi+\xi^{-1}$. Its Cartan/Schl\"afli matrix is $M=\begin{pmatrix} 2 & 2-\ga^2 \\ 2-\ga^2 & 2 \end{pmatrix}$ and $|\La_w/\La_r|=p$.
\end{ex}

\begin{ex}
Consider the root system of type $B_2$. Since it is crystallographic, its coefficient ring is $\mathcal{O}=\Z$. Its root lattice is the $\Z$-span of the simple roots $\Pi=\{\al_1, \al_2\}$. The corresponding Cartan matrix is $M=\begin{pmatrix} 2 & -2 \\ -1 & 2 \end{pmatrix}$ and $|\La_w/\La_r|=2$. Observe that normalized root system of type $B_2$ gives the dihedral system $I_2(4)$.
\end{ex}


\subsection{Formal group laws} \label{int:fgrs}
Let $F$ be a commutative one dimensional {\it formal group law} (FGL) over a commutative unital ring $R$ (the coefficient ring of $F$).  By definition (see \cite[\S1.1]{LM07}), $F$ is a formal power series in two variables $F(x,y)\in R[[x,y]]$ which is commutative $F(x,y)=F(y,x)$, associative $F(F(x,y),z)=F(x,F(y,z))$ and has neutral element $F(x,0)=x$. Using these properties one can express $F$ as
\[
F(x,y)=x+y+\sum_{i,j>0} a_{ij}x^i y^j,\;\text{ for some coefficients }a_{ij}\in R.
\]

The quotient of the polynomial ring in $a_{ij}$s modulo relations imposed by the commutativity ($a_{ij}=a_{ji}$) and associativity gives the Lazard ring $\mathbb{L}$. There is a universal FGL over $\mathbb{L}$ denoted $F_{\mathbb{L}}$ such that for any FGL $F$ over $R$ there is a unique map $\mathbb{L} \to R$ given by evaluating the coefficients $a_{ij}$ of $F_{\mathbb{L}}$.

\begin{exs}\label{ex:fgls} 
All examples of finite (polynomial) FGLs can be obtained by the base change (extending the coefficient ring) from the following cases (see \cite{La55} and \cite{LM07}):

The {\em additive} FGL $F(x,y)=x+y$ over $R=\Z$ with $F_{\mathbb{L}}\to F$ given by $a_{ij}\mapsto 0$ for all $i,j>0$.

The {\em connective} FGL $F(x,y)=x+y-\be xy$ over $R=\Z[\be]$, $\be\neq 0$, where $F_{\mathbb{L}}\to F$ is given by $a_{11}\mapsto -\be$, $a_{ij}\mapsto 0$ for all $i+j>2$. 

If $\be$ is invertible and $R=\Z[\be^{\pm 1}]$, then the respective $F$ is called {\em graded multiplicative}. If $\be=1$, the FGL $F$ over $R=\Z$ is called {\em multiplicative}.
\end{exs}


\subsection{Formal group rings}\label{int:fgring}
Given a crystallographic root system $\Phi$, we fix 
\begin{itemize}
\item
a free intermediate $\Z$-module (lattice) $\La$, so that $\La_r \subset \La \subset \La_w$, and
\item
a FGL $F$ over a ring $R$.
\end{itemize}
Following \cite[\S2]{CPZ} we define the {\em formal group ring} (FGR) as follows:

We first consider the polynomial ring $R[x_\la]_{\la\in \La}$ in variables corresponding to elements of $\La$. We then take  its completion $R[[x_\la]]_{\la\in \La}$ with respect to the kernel $I$ of the augmentation map $\ep\colon x_\la\mapsto 0$. Finally, we set
\[
R[[\La]]_F:=R[[x_\la]]_{\la\in \La}/\overline{J},
\]
where $\overline{J}$ denotes the closure (in the $I$-adic topology) of the ideal of relations $x_0=0$ and $x_{\la+\mu}=F(x_\la,x_\mu)$, $\forall \la,\mu\in \La$. 

If $F$ is finite (see~\ref{ex:fgls}), we skip the completion step and set the FGR to be
\[
R[\La]_F:=R[x_{\la}]_{\la \in \La}/J.
\] 

From this point on we denote by $S$ either the FGR $R[\La]_F$ for a finite $F$ or the FGR $R[[\La]]_F$ for an infinite $F$. 

\begin{exs}
Depending on a choice of finite FGL (see~\ref{ex:fgls}) we obtain the following examples of FGRs:

\medskip

{\em Additive:} $S=Sym_\Z(\La)$ is the symmetric algebra over $R=\Z$. After fixing a basis $\{x_1,\ldots,x_n\}$ of $\La$ we may identify $S$ with the polynomial ring $S=\Z[x_1,\ldots,x_n]$ (here $x_\la$ corresponds to $\la$).

\medskip

{\em Connective:} $S=R[x_\la]_{\la \in \La}/(x_0,x_{\la+\mu}-x_{\la}-x_\mu+\be x_\la x_\mu)$, where $R=\Z[\be]$. The ring $S$ can be turned into a graded ring over $\Z$ by setting $\deg x_\la=1$ and $\deg \be =-1$. Let $I=\langle (1-e^\la)\rangle_\la$ be the augmentation ideal of the integral group ring $\Z[\La]$. Then the degree 0 component of $S$ can be identified with the Rees ring $\Z[\La][I]$ via $1-e^\la \mapsto \be x_{-\la}$. 

\medskip

{\em (Graded) multiplicative:} $S=\Z[\La]\otimes_\Z R$ is the group ring over Laurent polynomials $R=\Z[\be^{\pm 1}]$ (here $e^\la\mapsto 1-\be x_{-\la}$). This is a graded ring as in the connective case. If $\be=1$ we obtain the usual integral group ring $S=\Z[\La]$.
\end{exs}

In the case of an additive FGL we extend the definition of a FGR $S$ to any non-crystallographic root system $\Phi$  as follows:

Consider the additive FGL over the coefficient ring $R:=\mathcal{O}$ of $\Phi$ and set $S$ to be the symmetric algebra over $\mathcal{O}$ of the respective free $\mathcal{O}$-submodule $\La$, that is
\[
S=Sym_{\mathcal{O}}(\La)=\mathcal{O}[x_{1},\ldots, x_{n}],
\]
where $\La_r\subset \La\subset \La_w$ and $\{x_1,\ldots,x_n\}$ is an $\mathcal{O}$-basis of $\La$.


\subsection{Push-pull elements}\label{int:pushpull} 
We now recall basic notation from \cite{CZZ1}. Let $W$ be the respective (finite) Coxeter/Weyl group generated by reflections 
\[
s_\al\colon \la \mapsto \la - \al^\vee(\la)\al,\qquad \al\in \Phi_+,\;\la\in \La.
\] 
It acts on the FGR $S$ through $\La$ by reflections. We denote by $S_W$ the twisted group algebra that is the left free $S$-module $S\otimes_{R} R[W]$ with basis 
$\{\de_w\}_{w\in W}$ and multiplication induced by the relation 
\[
\de_w f = w(f)\de_w,\quad f\in S.
\]
To proceed to the next step we need to assume that each $x_{\al}$, $\al\in \Phi_+$, is regular in $S$, that is a non-zero element which is not a zero-divisor. Observe that modulo the $W$-action this is equivalent to assume that  $x_\al$ is regular,  where $\al\in \Pi$.

\begin{rem}
According to \cite[Lemma~2.2]{CZZ1} $x_{\al}$ is always regular except when $\al$ is a long root in the irreducible root system of simply connected type $\mathrm{C}$. In the latter case it is regular if and only if  $2\cdot_F x$ is regular in $R[[x]]$ which always holds if $2$ is regular in $R$.
\end{rem}

Let $Q$ be the localization $S[\tfrac{1}{x_{\al}}]_{\al\in \Phi^+}$ of $S$ at all $x_\al$'s, and let $Q_W$ denote the respective twisted group algebra.
For each root $\al$ we define the {\em push-pull element} in $Q_W$ by
\[
Y_\al:=\tfrac{1}{x_{-\al}}+\tfrac{1}{x_\al}\de_{s_\al}.
\]

The twisted group algebra $Q_W$ acts on $Q$ via $q\de_w (f)=q w(f)$, $q,f\in Q$. Under this action the push-pull elements correspond to the classical push-pull operators
\[
Y_\al(f)=\tfrac{f}{x_{-\al}}+\tfrac{s_\al(f)}{x_\al},\quad f\in S
\]
which satisfy the quadratic relation:
\[
Y_\al^2=\kappa_\al Y_\al,\quad\text{ where }\kappa_\al=\tfrac{1}{x_\al}+\tfrac{1}{x_{-\al}}\in S.
\]
We may also consider {\em divided-difference elements} $X_\al:=\tfrac{1}{x_\al}-\tfrac{1}{x_\al}\de_{s_\al}$. These elements correspond to the usual divided difference operators  $X_\al(f)=\tfrac{f-s_\al(f)}{x_\al}$, $f\in S$, and behave similar to $Y_i$s. Observe that in the additive case $X_\al=-Y_\al$.

\begin{ex} In the finite case we have
\[
x_{-\al}=\begin{cases}
-x_\al, \text{ in the additive case, and} \\
\tfrac{-x_\al}{1-\be x_\al}, \text{ in the connective case,}
\end{cases}
\]
so we obtain
\[
\kappa_\al=
\begin{cases}
0, \text{ in the additive case, and} \\
\be, \text{ in the connective case.}
\end{cases}
\]
\end{ex}


\subsection{Reduced words}
Given a simple root $\al_i$ we denote $s_i=s_{\al_i}$, $\de_i=\de_{s_i}$, $x_i=x_{\al_i}$, $x_{-i}=x_{-\al_i}$, $\kappa_i=\kappa_{\al_i}$ and $Y_i=Y_{\al_i}$, $Y_{-i}=Y_{-\al_i}$. The celebrated result by Bressler-Evens \cite{BE90} says that the push-pull elements satisfy usual braid relations in $W$ if and only if $F$ is finite (or equivalently, $\kappa_i \in R$ is a constant). Therefore, in the finite case we may set $Y_w:=Y_{i_1}\ldots Y_{i_l}$ for any reduced expression $w=s_{i_1}\ldots s_{i_l} \in W$ of length $l=\ell(w)$. In the infinite case, for each $w\in W$ and its reduced expression $w=s_{i_1}\ldots s_{i_l}$ we denote the respective product of push-pull elements by $Y_{I_w}$, where $I_w=(i_1,\ldots, i_l)$.

By definition, for push-pull elements we have
\[
\de_wY_\al \de_{w^{-1}}=Y_{w(\al)},\quad \al\in \Phi,\; w\in W,
\]
and for the respective push-pull operators we have
\[
{}^wY_{\al}(f)=wY_\al w^{-1}(f)=Y_{w(\al)}(f),\quad w\in W.
\]

Let $w_0$ be the element of maximal length in $W$. Observe that $-w_0$ acts by an automorphism (possibly identical) on the set of simple roots $\Pi$. By means of this automorphism a reduced word for $w$ is mapped to a reduced word for $w_0ww_0$. So for $I_w=(i_1,\ldots, i_l)$ we have 
\[
{}^{w_0} Y_{I_w}=Y_{w_0(\al_{i_1})}\ldots Y_{w_0(\al_{i_l})}=Y_{w_0(I_{w})}.
\]
Observe that in the finite case $Y_{w_0(I_{w})}=Y_{-I_{w_0ww_0}}$, where $-I_w=(-i_1,\ldots,-i_l)$.

\begin{ex} 
In the additive case, we have $Y_{-i}=-Y_i=X_i=\tfrac{1-\de_i}{\al_i}$.

In the connective case, we have
\[
Y_{-i}=\tfrac{1}{x_i}+\tfrac{1}{x_{-i}}\de_i=\tfrac{1}{x_i}+\tfrac{1-\be x_i}{-x_i}\de_i=X_i+\be\de_i.
\]
\end{ex}


\subsection{Structure algebra}\label{int:struct}
By a {\it structure algebra} $\SA$ of the Bruhat poset $(W,\le)$ we call a $S$-submodule of the free module $\oplus_{v\in W}S$ defined by
\[
\SA=\Big\{(z_v)_v\in \bigoplus_{v\in W} S\mid z_{s_\al w}-z_w \in x_\al S,\; \forall \; w\le s_\al w,\, \al\in \Phi_+ \Big\}
\] 
with the coordinate-wise multiplication. 

The structure algebra $\SA$ plays an important role in the theory of sheaves on moment graphs. Namely, consider the Bruhat poset $(W,\le)$ of \cite[II.5.9]{Hu90}. The data 
\[
\text{Vertices } V:=W,\quad \text{Labelled Edges } E:=\big\{w\stackrel{\al}\longrightarrow s_\al w\mid  w\le s_\al w,\, \al\in \Phi_+,\, w\in W\big\}
\]
define a labelled directed graph called the {\it Bruhat moment graph} and denoted $\Gg$. By definition such a graph does not have multiple edges, self-loops and directed cycles. The algebra $\SA$ then gives the ring of global sections of the structure sheaf on $\Gg$ that is
\[
\Gamma(\Gg;S)=\{(z_v)_{v\in V} \mid z_{y}-z_x \in x_\al S,\; x\stackrel{\al}\longrightarrow y \in E\}.
\]

\begin{ex}\label{ex:momentd}
In the dihedral case $I_2(p)$, the Coxeter group $W$ (symmetry group of a regular $2p$-gon) is generated by simple reflections $s_1, s_2$ corresponding to simple roots $\al_1=1, \al_2=\xi^{p-1}$. The respective Bruhat moment graph $\Gg$ is given by a regular $2p$-gon with two opposite vertices $e$ and $w_0$ (element of maximal length) and edges directed from $e$ to $w_0$ such that each $(2i+1)$-diagonal, $i=0\ldots \tfrac{p-1}{2}$ is labelled either by $\xi^i$ or $\xi^{p-1-i}$ depending on the direction. The respective algebra $\SA$ of global section is then a submodule of the free $S$-module of rank $2p$.
\end{ex}

On the other hand, the structure algebra is closely related to affine Demazure algebras and Hecke-type rings of \cite{HMSZ} and of \cite{CZZ1}. Indeed, following \cite{CZZ1} consider the $R$-subalgebra $\mathbb{D}$ of the twisted group algebra $Q_W$ generated by all the push-pull elements $Y_i$'s and elements of $S$. It is called the {\em affine Demazure algebra}. Observe that the products of push-pull elements $\{Y_{I_w}\}_{w\in W}$ form an $S$-basis of $\mathbb{D}$ (as a left $S$-module).

\begin{ex}
In the additive case $\mathbb{D}$ is the nil (affine) Hecke ring of \cite{KK86} and in the multiplicative case $\mathbb{D}$ coincides with the $0$-affine Hecke ring of \cite{KK90} respectively. 
\end{ex}

According to one of the main results of \cite{CZZ2} (Theorem 9.2 loc.cit.) if $S$ satisfies the conclusion of \cite[Lemma~2.7]{CZZ2}, that is 
\begin{equation}\tag{*}
\text{for any two distinct positive roots $\al$ and $\al'$ and $f\in S$},\;
x_\al\mid x_{\al'} f \Longrightarrow x_\al \mid f,
\end{equation}
then the structure algebra $\SA$ is isomorphic to the $S$-linear dual of $\mathbb{D}$. Namely, $\SA$ coincides with the image of the inclusion (the moment map)
\[
Hom_S (\mathbb{D},S)\longrightarrow Hom_S (S_W,S)=\oplus_{v\in W} S
\]
induced by $S_W\hookrightarrow \mathbb{D}$, $\de_i\mapsto (x_i Y_i - \tfrac{x_i}{x_{-i}})$. In what follows we will extensively use this fact, hence, assuming $S$ satisfies (*).

\begin{rem} 
Observe that in the crystallographic case the assumption~(*) on $S$ follows from regularity of certain formal integers in $R$ given in \cite[Table~1]{CZZ2}). Also this assumption immediately implies that each $x_\al$ is regular in $S$. In non-crystallographic cases $I_2(p)$, i.e., in the additive case, it follows from regularity of the prime $p \in \mathcal{O}$.
\end{rem}


\subsection{The Hecke action and the Weyl action}
There is an $S$-linear action of the twisted group algebra $Q_W$ on the structure algebra $\SA$ given by
\[
q\de_w \bullet (z_v)_v=(z_v')_v,\quad z_v'=v(q)z_{vw}.
\]
It is called the {\it Hecke action}. It satisfies $X\bullet Y \bullet z=(XY)\bullet z$, $X,Y\in Q_W$, $z\in \SA$. Under this action the push-pull elements act by
\[
Y_i\bullet (z_v)_v =(\tfrac{1}{x_{-i}}+\tfrac{1}{x_i}\de_i)\bullet (z_v)_v =
(\tfrac{z_v}{x_{vs_i(\al_i)}}+\tfrac{z_{vs_i}}{x_{v(\al_i)}})_v.
\]
Observe that in $(z'_v)_v=Y_i\bullet (z_v)_v$ coordinates $z_w'$ and $z_v'$ coincide if $w=vs_i$. 

\begin{ex} 
For example, we have
\[
Y_i\bullet (z_v)_v=
\begin{cases}
(\tfrac{(z_{vs_i}-z_v)}{v(\al_i)})_v, \text{ in the additive case,}\\
(\tfrac{(z_{vs_i}-z_v)}{x_{v(\al_i)}}+\be z_v)_v, \text{ in the connective case.}
\end{cases}.
\]
\end{ex}

The twisted group algebra $Q_W$ also acts on $\SA$ via the {\it Weyl action} by
\[
q\de_w \odot (z_v)_v=(z_v')_v, \; \text{ where } z_v'=qw(z_{w^{-1}v}).
\]
By definition, we have $s\odot z=sz$, $z\in \SA$, $s\in S$. Moreover, it commutes with the Hecke action, i.e.,
\[
X\bullet Y\odot z=Y\odot X\bullet z, \quad X,Y\in Q_W,\; z\in \SA.
\]
As for the push-pull elements, we have
\[
Y_i\odot (z_v)_v=(\tfrac{1}{x_{-i}}+\tfrac{1}{x_i}\de_i)\odot (z_v)_v =(\tfrac{z_v}{x_{-i}}+\tfrac{s_i(z_{s_iv})}{x_{i}})_{v}.
\]
\begin{ex}
For example, we have 
\[
Y_i\odot (z_v)_v=
\begin{cases}
(\tfrac{(s_i(z_{s_iv})-z_v)}{\al_i})_v, \text{ in the additive case,}\\
(\tfrac{(s_i(z_{s_iv})-z_v)}{x_i}+\be z_v)_v, \text{ in the connective case.}
\end{cases}.
\]
\end{ex}

We define the {\it characteristic map} $c\colon S\to \SA$ by 
\[
c(f)=f\bullet \mathbf{1}=(w(f))_w,\;\text{ where }\mathbf{1}=(1)_w.
\]
We also define the {\it Borel map} $\rho\colon S\otimes_{S^W} S \to \SA$ by 
\[
\rho(f_1\otimes f_2)=f_1\odot f_2\bullet \mathbf{1}=f_2\bullet f_1\odot \mathbf{1}=f_1c(f_2)=(f_1w(f_2))_w.
\]

Then for each $X\in Q_W$ we have two commutative diagrams
\[
\xymatrix{
S\otimes_{S^W} S \ar[r]^-\rho \ar[d]_{(X\cdot - )\otimes \id} & \SA\ar[d]^{X\odot -} \\
S\otimes_{S^W} \ar[r]^-\rho S & \SA
}\qquad 
\xymatrix{
S\otimes_{S^W} S \ar[r]^-\rho \ar[d]_{\id\otimes (X\cdot - )} & \SA\ar[d]^{X\bullet -} \\
S\otimes_{S^W} \ar[r]^-\rho S & \SA
}
\]
In particular, we have
\[
c(Y_{I_w}(f))=Y_{I_w}\bullet c(f).
\]


\section{Intersection multiplicities of twisted Schubert classes}\label{sec:intmult}

It is well-known (see e.g. \cite{CZZ1}) that the structure algebra $\SA$ (the respective cohomology $\hh_T(G/B)$) has a natural basis given by Schubert classes. In~\ref{subs:schub} we recall the definition, basic properties and examples of such classes. In geometric terms, these are obtained by applying the push-pull operators $Y_i$ to the classes of $T$-fixed points $[v]$. Observe that the natural $W$-action (the $\odot$-action) permutes such $[v]$'s giving rise to the so-called twisted Schubert classes.

As the next step in~\ref{subs:pushf} we study the action by the push-forward element $\pi\bullet\text{-}$ which is the algebraic counterpart of the classical equivariant push-forward map $\hh_T(G/B) \to \hh_T(pt)$. Observe that in the finite case this action coincides with the action by the operator $Y_{w_0}$ (the composite of push-pull operators corresponding to the element of maximal length).

Taking  the Poincar\'e duals of the Schubert classes gives the so-called dual basis. In the last subsection~\ref{subs:intmu} we study the transition matrices $C$ and $D$ between these two bases and their $u$-twisted versions for $u\in W$. Geometrically speaking, the entries of these matrices are given by (equivariant) intersection multiplicities. These had been intensively studied in the additive and multiplicative cases.

In the additive case (the case of a usual cohomology), these matrices can be viewed as a deformation from the Littlewood-Richardson coefficients (case $u=e$) to the complimentary diagonal matrix (case $u=w_0$) which reflects the fact that Poincar\'e duals of Schubert classes are again Schubert classes. In the connective case, the case $u=w_0$ gives the transition matrix between the Schubert basis and Brion's basis of ideal sheaves in $K$-theory \cite{Br02}. In Example~\ref{ex:connm} we show how to compute such matrix $C_{w_0}$.  

In the case of a general FGL, the definition and basic properties of such matrices (intersection multiplicities) have never been discussed or studied directly, though can be easily deduced from the results of \cite{CZZ2}.  


\subsection{Classes of fixed points and twisted Schubert classes}\label{subs:schub}
Consider the element $x_\Pi=\prod_{\al\in \Phi^-}x_\al \in S$ (the product over all negative roots). Set $[e]=(z_v)_v \in \SA$, where $z_{e}=x_\Pi$ and all other coordinates $z_v$ are zeros. We call $[e]$ the class of a fixed point at $e$. We then set $[w]=\de_w\odot [e]$ and call it the class of a fixed point at $w$. It has $z_w=w(x_\Pi)$ and all other coordinates are zeros.

\begin{dfn}
We define the Schubert classes at $e$ (we call it the usual Schubert classes) via
\[
\zeta_{I_w}:=Y_{I_w^{-1}} \bullet [e],\quad w\in W.
\]
Similarly, we define the Schubert classes at $w_0$ (we call it special/opposite Schubert classes) by
\[
\sigma_{I_w}:=Y_{I_w^{-1}} \bullet [w_0],\quad w\in W.
\]
\end{dfn}

\begin{rem}
Suppose we are in the additive case, so $S$ is a polynomial ring. Given $w\in W$, an element $\theta^w=(z_v)_v\in \SA$ is called special if there is some fixed $d$ so that each $z_v$ is a homogeneous polynomial of degree $d$ and $\theta^w$ satisfies the following support condition:
\[
z_v\neq 0\quad \Longleftrightarrow\quad w\le v.
\]
It is known that for any special element $\theta^w$ the coordinate $z_w$ is divisible by $w(x_\Pi)$. Moreover, there are unique special elements where this coordinate is precisely $w(x_\Pi)$. These elements are given precisely by the special Schubert classes $\sigma_w$.
\end{rem}

More generally, we define the Schubert classes at any $v\in W$ (and call it twisted Schubert classes) by 
\[
\zeta(v)_{I_w}:=Y_{I_w^{-1}}\bullet [v],\quad w\in W.
\]

By the very definition we have $\zeta(e)_{I_w}=\zeta_{I_w}$ and $\zeta(w_0)_{I_w}=\sigma_{I_w}$. In all the finite FGL cases we will simply write $\zeta(v)_w$ instead of $\zeta(v)_{I_w}$.

Observe that there is an anti-involution on $\mathbb{D}$ given by (see \cite[\S3.1.(6)]{LZZ})
\[
\imath\colon q\de_w \mapsto \de_{w^{-1}}q \tfrac{w(x_\Pi)}{x_\Pi}
\] 
which switches $\odot$ and $\bullet$ if applied to the class $[e]$ and such that $\imath(Y_{I_w})=Y_{I_w^{-1}}$. Applying it we obtain
\[
\zeta(e)_{I_w}=Y_{I_w^{-1}} \bullet [e] = \imath(Y_{I_w^{-1}})\odot [e]=Y_{I_w}\odot [e],
\]
and, more generally,
\[
\zeta(v)_{I_w}=Y_{I_w^{-1}}\bullet [v]=Y_{I_w^{-1}} \bullet \de_v \odot [e]=\de_v\odot \zeta_{I_w}.
\]

\begin{ex} 
In the $SL_2$-case we have $\Pi=\{\al_1\}$, $W=\{e, s_1=w_0\}$, $x_\Pi=x_{-1}$. The only push-pull operator $Y_1=\tfrac{1}{x_{-1}}+\tfrac{1}{x_1}\de_1$ acts as follows:
\[
Y_1\bullet (z_e,z_{1})=(\tfrac{z_e}{x_{-1}}+\tfrac{z_1}{x_1},\tfrac{z_e}{x_{-1}}+\tfrac{z_1}{x_1}).
\]
In particular, we obtain
\[
Y_1\bullet [e] =Y_1\bullet (x_{-1},0)=\mathbf{1}, \; \text{ and }\; 
Y_1\bullet [w_0] =Y_1 \bullet (0, x_1)=\mathbf{1}.
\]
\end{ex}

\begin{ex}
In the $SL_3$-case, we have $\Pi=\{\al_1,\al_2\}$, $x_\Pi=x_{-\al_1}x_{-\al_2}x_{-\al_1-\al_2}$, and $W=\{e, s_1,s_2,s_2s_1,s_1s_2,s_1s_2s_1=s_2s_1s_2=w_0\}$. There are two push-pull operators $Y_i=\tfrac{1}{x_{-i}}+\tfrac{1}{x_i}\de_i$, $i=1,2$ which act as follows ($i\neq j$):
\[
Y_i\bullet (z_e,z_{i},z_j,z_{ji},z_{ij},z_{w_0})=(z_e',z_i'=z_e',z_j',z_{ji}'=z_j',z_{ij}',z_{w_0}'=z_{ij}'),
\]
where $z_e'=\tfrac{z_e}{x_{-i}}+\tfrac{z_i}{x_i}$, $z_j'=\tfrac{z_j}{x_{-i-j}}+\tfrac{z_{ji}}{x_{i+j}}$, and $z_{ij}'=\tfrac{z_{ij}}{x_{-j}}+\tfrac{z_{w_0}}{x_j}$. Iterating this formula, we obtain
\begin{align*}
Y_i\bullet [e] &=(x_{-j}x_{-i-j},x_{-j}x_{-i-j},0,0,0,0), \\
Y_jY_i\bullet [e] &=(x_{-i-j},x_{-j},x_{-i-j},x_{-j},0,0), \\
Y_iY_jY_i\bullet [e] &=(\tfrac{x_{-i-j}}{x_{-i}}+\tfrac{x_{-j}}{x_i}, \tfrac{x_{-i-j}}{x_{-i}}+\tfrac{x_{-j}}{x_i},1,1,1,1).
\end{align*}
Applying $\de_{w_0}\odot -$, we obtain
\begin{align*}
Y_i\bullet [w_0] &=(0,0,0,0,x_ix_{i+j},x_ix_{i+j}), \\
Y_jY_i\bullet [w_0]&=(0, 0,x_i,x_{i+j},x_i ,x_{i+j}), \\
Y_iY_jY_i\bullet [w_0] &=(1,1,1,1,\tfrac{x_i}{x_{-j}}+\tfrac{x_{i+j}}{x_j},\tfrac{x_i}{x_{-j}}+\tfrac{x_{i+j}}{x_j}).
\end{align*}
Observe that $\tfrac{x_i}{x_{-j}}+\tfrac{x_{i+j}}{x_j}=1$ $\Longleftrightarrow$ the braid relation of type $A$ holds for $Y_i$s $\Longleftrightarrow$ $F$ is finite (see \cite{HMSZ}). So, in the additive and connective cases we obtain 
\[
Y_{w_0}\bullet [e]=Y_{w_0}\bullet [w_0]=\mathbf{1}.
\]
\end{ex}


\subsection{The push-forward element}\label{subs:pushf}
Following \cite[Def.~5.3]{CZZ2} we define the push-forward element $\pi\in Q_W$ by
\[
\pi=(\sum_{w\in W} \de_w)\tfrac{1}{x_\Pi}=\sum_{w\in W} \tfrac{1}{w(x_\Pi)}\de_w.
\]
It correspond to the push-forward operator $\pi(f)=\sum_w \tfrac{w(f)}{w(x_\Pi)}$, $f\in S$. Observe that it commutes with the characteristic map 
\[
\pi\bullet c(f)=c(\pi(f)),\quad  f\in S.
\]

The push-forward element $\pi$ acts on the structure algebra $\SA$ as
\[
\pi \bullet (z_v)_v =(\sum_w \tfrac{z_{vw}}{vw(x_\Pi)})_v,\;\text{ and }\; \pi\odot (z_v)_v=(\sum_w w(\tfrac{z_{w^{-1}v}}{x_\Pi}))_v.
\]

For the fundamental class we get 
\[
\pi\bullet \mathbf{1}=(\sum_w \tfrac{1}{w(x_\Pi)})\cdot \mathbf{1},\quad\text{ a multiple of }\mathbf{1}.
\]
For the class of a point $[v]=\de_v\odot [e]$ we obtain
\[
\pi\bullet (\de_v\odot [e])=\de_v\odot \pi\bullet [e] =\de_v\odot\mathbf{1}=\mathbf{1}.
\]

According to \cite[Lemma~7.1]{CZZ2} we have
\[
\pi \bullet \Big(\big(Y_i\bullet z \big)\cdot z'\Big)=\pi\bullet \Big(z \cdot \big(Y_i\bullet z'\big)\Big),\qquad z,z'\in \SA. 
\]
Combining it with the formula $\pi \bullet ([v] \cdot (z_w)_w)=z_{v}\mathbf{1}$, we obtain the following

\begin{lem}
For any $z\in \SA$ and $v,w\in W$ we have
\[
\pi\bullet (\zeta(v)_{I_w} \cdot z)=\pi\bullet \big([v]\cdot (Y_{I_w}\bullet z)\big)=(Y_{I_w}\bullet z)_v\mathbf{1}.
\]
\end{lem}


\subsection{The dual basis and intersection multiplicities}\label{subs:intmu}
According to \cite[\S15]{CZZ2} there is a perfect pairing $\SA\otimes_S \SA \mapsto S\mathbf{1}$, $(z,z')\mapsto \pi\bullet (z\cdot z')$ so that
\[
\pi\bullet (\zeta_{I_w}\cdot \zeta^\vee_{I_u})=\de^{Kr}_{w,u}\mathbf{1},
\]
where 
\[
\zeta^\vee_{I_w}=Y_{I_w}^* \in \SA=Hom_S(\mathbb{D},S),
\]
is the $S$-linear dual of $Y_{I_w}$ defined by $Y_{I_w}^*(Y_{I_u})=\de^{Kr}_{w,u}$. We call the classes $\zeta^\vee_{I_w}$ the Poincar\'e duals of $\zeta_{I_w}$.

\begin{lem} 
For twisted Schubert classes we have
\[
\zeta(v)_{I_w}^\vee =\de_v\odot \zeta_{I_w}^\vee.
\]
In particular, the classes $\{\zeta(v)_{I_w}^\vee\}$ form an $S$-basis of $\SA$.
\end{lem}

\begin{proof}
Indeed, since $\bullet$ and $\odot$ actions commute, we obtain
\[
\pi\bullet (\zeta(v)_{I_w}\cdot (\de_{v}\odot \zeta_{I_w}^\vee))=\pi\bullet (\de_v\odot (\zeta_{I_w}\cdot \zeta^\vee_{I_u}))=\de_v\odot \pi\bullet (\zeta_{I_w}\cdot \zeta^\vee_{I_w})=\de^{Kr}_{w,u}\mathbf{1}.\qedhere
\]
\end{proof}

\begin{dfn}
We define the (equivariant) intersection multiplicities $c_{w,v}^u\in S$ and $d_{w,v}^u\in S$, $w,v,u\in W$ by
\[
\pi\bullet (\zeta_{I_w} \cdot \zeta(u)_{I_v})=c_{w,v}^u\mathbf{1},\; \text{ and } \pi\bullet (\zeta_{I_w}^\vee \cdot \zeta(u)_{I_v}^\vee)=d_{w,v}^u\mathbf{1}.
\]
\end{dfn}

Applying $\de_{u^{-1}}\odot -$ to the both sides, we obtain
\[
u(c_{w,v}^{u^{-1}})=c_{v,w}^u\;\text{ and }u(d_{w,v}^{u^{-1}})=d_{v,w}^u.
\]
Therefore, if $C_u=(c_{w,v}^u)_{w,v}$ and $D_u=(d_{w,v}^u)_{w,v}$ denote the respective matrices of coefficients, then
\[
u(C_{u^{-1}})=C_u^t,\; \text{ and }\; u(D_{u^{-1}})=D_u^t.
\]

Observe also that by duality we have
\[
\zeta(u)_{I_v}=\sum_w c_{w,v}^u \zeta^\vee_{I_w}, \text{ and }\zeta(u)^\vee_{I_v}=\sum_w d_{w,v}^u \zeta_{I_w}.
\]
Therefore,
\[
\zeta_{I_v}^\vee=\de_{u} \odot \zeta(u^{-1})^\vee_{I_v}=\de_u\odot \sum_w d_{w,v}^{u^{-1}} \zeta_{I_w}=
\sum_w u(d_{w,v}^{u^{-1}})\zeta(u)_{I_w}
\]
which implies $C_u^{-1}=D_u^t$.

Recall that we have
\[
\pi\bullet (\zeta_{I_w} \cdot \zeta(u)_{I_v})=\pi\bullet ([e] \cdot Y_{I_w}Y_{I_v^{-1}}\bullet [u])=(Y_{I_w}Y_{I_v^{-1}}\bullet [u])_e,
\]
which immediately gives (by looking at $e$-coordinates of $Y_{I_w}Y_{I_v^{-1}}\bullet [u]$)
\[
c_{w,v}^{u}=0 \text{ if } l(w)+l(v)<l(u).
\]
Similarly, we get
\[
d_{w,v}^{u}=0 \text{ if }l(w_0w)+l(w_0v)<l(u).
\]

\begin{ex} 
In the $SL_2$-case, we have $\pi=Y_1$. For the Schubert classes $\zeta_e=[e]$ and $\zeta_{w_0}=\mathbf{1}$ we compute their duals as follows. Let $(z_e,z_1)=[e]^\vee$. Then we have
\[
Y_1\bullet ([e]\cdot (z_e,z_1))=Y_1\bullet (x_{-1}z_e,0)=(z_e,z_e)=(1,1),
\]
which implies $z_e=1$. At the same time
\[
Y_1\bullet(\mathbf{1} \cdot (1,z_1))=Y_1\bullet (1,z_1)=(\tfrac{1}{x_{-1}}+\tfrac{z_1}{x_1},\tfrac{1}{x_{-1}}+\tfrac{z_1}{x_1})=(0,0),
\]
which implies $z_1=\tfrac{-x_1}{x_{-1}}$. Hence, $[e]^\vee=(1,\tfrac{-x_1}{x_{-1}})$ and, therefore,
\[
[e]^\vee=\kappa_1\tfrac{x_1}{x_{-1}}(x_{-1},0)+\tfrac{-x_1}{x_{-1}}(1,1)=\kappa_1\tfrac{x_1}{x_{-1}}[e]+\tfrac{-x_1}{x_{-1}}\mathbf{1}.
\]
Similarly, we find $\mathbf{1}^\vee=(0,x_1)$ and, therefore,
\[
\mathbf{1}^\vee=\tfrac{-x_1}{x_{-1}}(x_{-1},0)+x_1(1,1)=\tfrac{-x_1}{x_{-1}}[e]+x_1\mathbf{1}.
\]
These give 
$D_e=\begin{pmatrix} \kappa_1\tfrac{x_1}{x_{-1}} & \tfrac{-x_1}{x_{-1}} \\ \tfrac{-x_1}{x_{-1}} & x_1 \end{pmatrix}$
and $C_e=\begin{pmatrix} x_{-1} & 1 \\
1 & \kappa_1 \end{pmatrix}$ where $\kappa_1=\tfrac{1}{x_{-1}}+\tfrac{1}{x_1}$. Similarly, we get
$D_{w_0}=\begin{pmatrix} -\kappa_1 & 1 \\
1 & 0 \end{pmatrix}$ and $C_{w_0}=\begin{pmatrix} 0 & 1 \\ 
1 & \kappa_1 \end{pmatrix}$.
\end{ex}

\begin{ex}\label{ex:addc}
In the additive case, since $Y_i$'s generate the nil-Hecke ring ($Y_i^2=0$), $Y_{I_w}Y_{I_v^{-1}}=0$ and, hence, $c_{w,v}^u=0$ if $l(w)+l(v)>l(wv^{-1})$. In the case $l(w)+l(v)=l(wv^{-1})$ (or equivalently when the concatenation $I_wI_{v}^{-1}$ is reduced), we have $Y_{I_w}Y_{I_v^{-1}}=Y_{wv^{-1}}$ and
\[
c_{w,v}^u=(\zeta(u)_{vw^{-1}})_e=(\de_u\odot \zeta_{vw^{-1}})_e=u^{-1}(\zeta_{vw^{-1}})_u.
\]
In particular, for $u=w_0$ it implies that
\[
c_{w,v}^{w_0}=d_{w,v}^{w_0}=\de^{Kr}_{w_0w,v}\;\text{ for all }w,v\in W. 
\]
In other words, the special Schubert classes are Poincar\'e dual to the usual Schubert classes:
\[
\zeta_{w}^\vee=\sigma_{w_0w}=\de_{w_0}\odot \zeta_{w_0w}.
\]
\end{ex}

\begin{ex}\label{ex:connm}
In the connective case, $Y_i^2=\be Y_i$ in $\mathbb{D}$ (observe that for $\be=1$, $\mathbb{D}$ is the affine 0-Hecke algebra). Let $\overline{I}$ denote the reduced word in the 0-Hecke algebra. Then $Y_{I}=\be^{|I|-|\bar I|}Y_{\bar I}$ in $\mathbb{D}$ for any word $I$. Therefore, 
\[
c_{w,v}^u=\begin{cases}
\be^{|I_wI_{v}^{-1}|-|\overline{I_wI_{v^{-1}}}|}(Y_{\overline{I_w I_{v}^{-1}}}\bullet [u])_e,\; \text{ if } |\overline{I_wI_{v^{-1}}}|\ge \ell(u), \text{ and}\\
0, \text{ otherwise.}
\end{cases}
\]
In particular, $c_{w,v}^{w_0}=\be^{|I_wI_{v}^{-1}|-|\overline{I_wI_{v^{-1}}}|} \de^{^{Kr}}_{|\overline{I_wI_{v^{-1}}}|,l(w_0)}$. If $r=|\Phi_+|$, then 
\[
c_{w,v}^{w_0}=\begin{cases}
\be^{l(w)+l(v)-r}\; \text{ if } \overline{I_wI_{v^{-1}}}=w_0, \text{ and}\\
0, \text{ otherwise.}
\end{cases}
\]

To completely determine the matrix for $u=w_0$ it is hence sufficient to give a criterion for $\overline{I_wI_{v^{-1}}}=w_0$. This is the content of \cite[Lemma 4.3]{Fa05}, which tells the following: say that a subexpression $I'=(i_{j_1}, \ldots, i_{j_t})$ of a reduced expression $I=(i_1, \ldots, i_r)$ is good if $I'$ is reduced and for any $u,k$ such that $j_{k-1}<u<j_k$,  we have $l(s_{i_{j_1}}\ldots s_{i_{j_{k-1}}}s_u) > l(s_{i_{j_1}}\ldots s_{i_{j_{k-1}}})$; then $\overline{I_wI_{v^{-1}}}=w_0$ if and only if $v=w_0s_{i_{j_1}}\ldots s_{i_{j_t}}$ for some good sub-expression $(i_{j_1},\ldots i_{j_t})$ of a reduced expression $(i_1\ldots i_r)$ for $w$.

For example, in the $SL_2$ and $SL_3$ case we obtain, respectively, the matrices
\[
C_{w_0}=\left(
\begin{array}{cc}
0&1\\
1&\be
\end{array}
\right)\quad\text{ and }\quad
C_{w_0}=\left(
\begin{array}{cccccc}
0&0&0&0&0&1\\
0&0&0& 0 &1&\be\\
0&0&0&1&0 &\be\\
0& 0&1&\be&\be&\be^2\\
0&1&0&\be&\be&\be^2\\
1&\be&\be&\be^2&\be^2&\be^3
\end{array}
\right)
\]
(here the rows/columns in the $A_2$ matrix correspond to $(e,s_1,s_2, s_2s_1, s_1s_2,w_0)$).
\end{ex}


\section{The push-forward of a product}\label{sec:pushforw}

In the present section we prove and discuss different variations and applications of the push-forward of the product formula for structure algebras.

First, we generalize the classical Demazure formula \cite{De73} to the oriented equivariant setup (Lemma~\ref{lem:demform}). Observe that it was proven in the oriented non-equivariant setup in \cite{CPZ}. As an immediate consequence we establish the formula for the fundamental class (Corollary~\ref{cor:fundcl}).

Then we prove one of our key results -- the push-forward formula of a product for arbitrary structure algebras (Theorem~\ref{thm:pushfres}). The latter was established for ordinary cohomology under certain restrictions in \cite{Xi21}. In all other cases it seems to be new. 

As the next step, in~\ref{subs:dschub} we introduce double Schubert elements. In the additive and multiplicative cases it gives the well-known double Schubert and Grothendieck polynomials respectively. Using these double elements and the Borel map we prove the reconstruction formula (Corollary~\ref{cor:recform}) and establish the antipode of the Demazure formula (Corollary~\ref{cor:antip}). The latter leads to several interesting combinatorial identities involving double Schubert elements and Schubert classes. It is also used later to construct the twisted coproduct.


\subsection{The push-forward formula}
We start with the following equivariant generalization of the Demazure formula.

\begin{lem}\label{lem:demform} 
We have the following expression for the characteristic map for $s\in S$
\[
c(u^{-1}s)=\sum_w {}^uY_{I_w}(s) \zeta(u)^\vee_{I_w},\;\text{ where }{}^u Y_{I_w}=uY_{I_w}u^{-1}.
\]
\end{lem}

\begin{proof}
Indeed, since $\{\zeta(u)_{I_w}^\vee\}_{w\in W}$ is a basis, we can write $c(s)=\sum_w p_w \zeta(u)^\vee_{I_w}$ for some $p_w\in S$. Applying $\pi\bullet (- \cdot \zeta(u)_{I_v})$ to the both sides we obtain
\[
\pi\bullet (c(s)\cdot \zeta(u)_{I_v})=\sum_w p_w \pi\bullet (\zeta(u)_{I_w}^\vee \cdot \zeta(u)_{I_v})=p_{v}\mathbf{1}.
\]
As for the left side we get
\[
\pi\bullet (c(s)\cdot\zeta(u)_{I_v})=\pi\bullet ([u]\cdot (Y_{I_v}\bullet c(s)))=
\pi\bullet ([u]\cdot c(Y_{I_v}(s)))=uY_{I_v}(s)\mathbf{1}.
\]
Comparing the left and the right sides we obtain the desired formula.
\end{proof}

Observe also that each of the summands on the right hand side depend on a choice of reduced expression, but the left hand side doesn't.

\begin{cor}\label{cor:fundcl}
We have the following expressions for the fundamental class:
\[
\mathbf{1}=[u]^\vee + \sum_{w\neq e}  {}^uY_{I_w}(1) \zeta(u)^\vee_{I_w} = \zeta_{I_{w_0}}+\sum_{v\neq w_0} (\sum_w d_{v,w}^u {}^u Y_{I_w}(1))  \zeta_{I_v}.
\]
\end{cor}

\begin{proof} 
Substituting $s=1$ we get the first equality. As for the second, expressing the dual $\zeta(u)^\vee_{I_w}$ in terms of Schubert classes we obtain
$\mathbf{1}=\sum_{w,v} d_{v,w}^u {}^u Y_{I_w}(1)  \zeta_{I_v}$. The expression $\sum_w d_{v,w}^u {}^u Y_{I_w}(1) $ gives the coefficient at $\zeta_{I_v}$ and, hence,
does not depend on a choice of $u$. For instance, for $v=w_0$ taking $u=w_0$ we obtain the coefficient $1$ at $\zeta_{I_{w_0}}$. 
\end{proof}

\begin{ex}
In the additive case we have ${}^u Y_w(1)=\de^{Kr}_{w,e}$ and $\mathbf{1}=\zeta_{w_0}$ which implies
$d_{v,e}^u=\de^{Kr}_{v,w_0}$.
In the connective case we have ${}^u Y_w(1)=\be^{\ell(w)}$ and $\mathbf{1}=\zeta_{w_0}$ which gives the system of equations
\[
\sum_{w} d_{v,w}^u \be^{\ell(w)}=\de^{Kr}_{v,w_0}.
\]
\end{ex}

We are now ready to prove our first main result 

\begin{thm}[The push-forward formula for a product]\label{thm:pushfres}
The following identity holds in the structure algebra $\SA$:
\[
\pi\bullet (c(f) \cdot c(g))=\sum_{w,v} Y_{I_w}(f) \cdot uY_{I_v}(g)\cdot  (\pi\bullet (\zeta_{I_w}^\vee\cdot \zeta(u)_{I_v}^\vee)), \quad f,g\in S,\; u\in W.
\]
\end{thm}

\begin{proof}
In the structure algebra $\SA$ we have
\begin{align*}
\pi\bullet (c(f)\cdot c(g)) &=\pi\bullet (\sum_w Y_{I_w}(f)\zeta_{I_w}^\vee  \cdot c(g))\\
&=\sum_w Y_{I_w}(f) \cdot \pi\bullet (\de_{u}\odot \zeta(u^{-1})_{I_w}^\vee \cdot c(g))\\
&=\sum_w Y_{I_w}(f) \cdot \pi\bullet (\sum_v u(d_{v,w}^{u^{-1}})\zeta(u)_{I_v} \cdot c(g))\\
&=\sum_{w,v} u(d_{v,w}^{u^{-1}}) Y_{I_w}(f) \cdot \pi\bullet \big([u] \cdot  Y_{I_v}(c(g)))\big)\\
&=\sum_{w,v} d_{w,v}^u Y_{I_w}(f) \cdot \pi\bullet \big([u] \cdot  c(Y_{I_v}(g)))\big)\\
&=\big(\sum_{w,v} d_{w,v}^u Y_{I_w}(f)u(Y_{I_v}(g))\big)\cdot \mathbf{1}.\\
\end{align*}
Substituting the expression for $d_{w,v}^u$ we obtain the desired formula.
\end{proof}

\begin{cor}[The push-forward formula at a point] 
The following identity holds in the ring $S$:
\[
\pi(f \cdot u^{-1}g)=\sum_{w,v} d_{w,v}^u \cdot Y_{I_w}(f)\cdot {}^uY_{I_v}(g), \quad f,g\in S,\; u\in W
\]
where both sides of the equality are $W$-invariant.
\end{cor}

\begin{proof}
The left side of the identity of the theorem can be written as
\[
\pi\bullet (c(f)\cdot c(u^{-1}g))=c(\pi(f\cdot u^{-1}g))=(y(\pi(f\cdot u^{-1}g)))_y.
\]
The formula then follows from the identities
\[
y(\pi(f\cdot u^{-1}g))=\sum_{w,v} d_{w,v}^u Y_{I_w}(f)uY_{I_v}u^{-1}(g),\; \text{for each }y\in W.\qedhere
\]
\end{proof}

\begin{ex} 
In the additive case for $u=w_0$ following Example~\ref{ex:addc} we obtain
\[
\pi(f \cdot w_0 g)=\sum_{w,v} d_{w,v}^{w_0} Y_w(f){}^{w_0}Y_v(g) =\sum_w Y_{w}(f) X_{ww_0}(g)=\sum_w (-1)^{\ell(w)}X_{w}X_{ww_0}(g).
\]
This gives the top Leibniz rule of \cite{Xi21}.
\end{ex}


\subsection{Double Schubert elements.}\label{subs:dschub}
Consider the usual (non-equivariant) characteristic map
\[
\ep c\colon S \longrightarrow \SA\longrightarrow \SA/I\SA.
\]
Let $[pt]=\ep ([e])$ denote the class of a (non-equivariant) point. According to \cite{CPZ} there exists an element $s \in I^{|\Phi_+|}$ such that $\ep c(s)=\tau [pt]$, where $\tau\in R$ is the torsion index of the respective root datum. Observe that $\tau$ is always divisible by the torsion index of the respective root system, and by the order $|\La_w/\La|$. We refer to \cite{CPZ} and \cite{CZZ1} for the discussion concerning the possible values of $\tau$. We only mention the following

\begin{exs} 
For weight lattices and crystallographic root systems we can choose $s$ and $\tau$ such that $\tau\in \Z$ coincides with the torsion index of $\Phi$ in the additive case, and $\tau=1$ in the multiplicative case.
\end{exs}

From this point on we assume that $\tau$ is regular in $R$ and that $R$ has no $2$-torsion. By \cite{CZZ1} the Borel map $\rho$ then becomes an isomorphism after inverting $\tau$. Let $z$ be a preimage of $[e]$ under $\rho$ (after inverting $\tau$). Clearing denominators we may assume that \[z_{[e]}:=\mu z \in S\otimes_{S_W} S,\] where $\mu=\tau^{r}$ for some $r>1$. Observe that by the very definition $\rho(z_{[e]})=\mu [e]$ and $\rho(u^{(1)} z_{[e]})=\mu [u]$. We denote $z_{[u]}=u^{(1)} z_{[e]}$. 

\begin{dfn}
Define $\SC_{I_w}:=Y_{I_w^{-1}}^{(2)}(z_{[e]})\in S\otimes_{S^W} S$ (here $Y_{I_w^{-1}}^{(2)}$ denotes $Y_{I_w^{-1}}$ applied to the second factor). So $\SC_e=z_{[e]}$ and
\[
\rho(\SC_{I_w})=\rho (Y_{I_w^{-1}}^{(2)}(z_{[e]}))=\mu Y_{I_w^{-1}}\bullet [e]= \mu\zeta_{I_w}.
\]
We call $\SC_{I_w}$ double Schubert elements.
\end{dfn}

By the very definition we have  $\rho(u^{(1)}(\SC_{I_w}))=\mu \zeta(u)_{I_w}$. So $u^{(1)}(\SC_{I_w})$ denoted also by $\SC(u)_{I_w}$ corresponds to the twisted Schubert class and will be called the twisted double Schubert element.


\subsection{The reconstruction formula}
Applying the push-forward formula at a point to the second components of $z_{[v]}=\SC(v)_e$ and $1\otimes u^{-1}g$, $g\in S$, we obtain 
\begin{align*}
\pi^{(2)}(z_{[v]} \cdot (1\otimes u^{-1}g))&=\sum_{w,v} d_{w,v}^uY_{I_w}^{(2)}(z_{[v]})\cdot  (1\otimes {}^uY_{I_v}(g))\\
 &=\sum_{w} \SC(v)_{I_w^{-1}} \cdot  (1\otimes \sum_v d_{w,v}^u {}^u Y_{I_v}(g))\\
 &=\sum_w \SC(v)_{I_w^{-1}} \cdot  (1\otimes  {}^u \tilde Y_{I_w}(g)),
\end{align*}
where $ {}^u\tilde Y_{I_w}=\sum_v d_{w,v}^u {}^uY_{I_v}$ in $\mathbb{D}$. On the other side 
\[
\rho \pi^{(2)}(z_{[v]}\cdot (1\otimes u^{-1}g))= \pi\bullet (\mu [v]\cdot c(u^{-1}g))= \mu vu^{-1} g \mathbf{1}.
\]
Since $\rho(vu^{-1} g\otimes 1)=vu^{-1}g\mathbf{1}$, $\tau$ is regular, and $\rho$ is an isomorphism after inverting $\tau$ we obtain
\begin{cor}[Reconstruction formula]\label{cor:recform} 
For any $g\in S$ and $u,v\in W$ we have in $S\otimes_{S^W} S$
\[
\mu (vu^{-1}g \otimes 1) = \sum_w \SC(v)_{I_w^{-1}} \cdot  (1\otimes {}^u \tilde Y_{I_w}(g)).
\]
\end{cor}

Applying the Borel map $\rho$ to both sides of the reconstruction formula we obtain the following `antipode` of the Demazure formula:

\begin{cor}\label{cor:antip} 
In the structure algebra $\SA$ for any $g\in S$ and $u,v\in W$ we have 
\[vu^{-1}g\mathbf{1}=\sum_w c({}^u\tilde Y_{I_w}(g)) \zeta(v)_{I_w^{-1}}.\] 
\end{cor}

\begin{ex}
In the additive case for $u=v=w_0$ we get ${}^{w_0} \tilde Y_w={}^{w_0}Y_{w_0w}=X_{ww_0}$. So we obtain
\[
g\mathbf{1}=\sum_w c(X_{ww_0}(g))\sigma_{w^{-1}}.
\]
\end{ex}

\begin{ex}
Combining the Demazure formula and its antipode analogue for $u=v$ we obtain
\begin{align*}
g\mathbf{1}&=\sum_{w,x} Y_{I_x}{}^u \tilde{Y}_{I_w}(g)  \zeta_{I_x}^\vee \zeta(u)_{I_w^{-1}}  \\
&=\sum_{w,x} Y_{I_x} (\sum_y d_{w,y}^u {}^u Y_{I_y})(g) (\sum_z u(d_{z,x}^{u^{-1}})\zeta(u)_{I_z})\zeta(u)_{I_w^{-1}}  \\
&=\sum_{w,z} (\sum_{x,y} u(d_{z,x}^{u^{-1}}) Y_{I_x} d_{w,y}^u {}^u Y_{I_y})(g))\cdot \zeta(u)_{I_z}\zeta(u)_{I_w^{-1}}.
\end{align*}
In the additive case for $u=v=w_0$  we obtain the following combinatorial identities involving products of Schubert classes
\begin{align*}
g\mathbf{1}&=\sum_{w,x} Y_x {}^{w_0}Y_{w_0w}(g)\cdot \sigma_{w_0x}\sigma_{w^{-1}} \\
&=\sum_{x,y} (-1)^{l(w_0x )}X_{w_0x} X_{y^{-1}w_0}(g) \cdot \sigma_{x} \sigma_{y} \\
&=\sum_{{ \begin{array}{c} x,y \\ l(x)+l(y)=l(xy^{-1})\end{array}}} (-1)^{l(w_0x )} X_{w_0xy^{-1}w_0}(g) \cdot \sigma_{x} \sigma_{y}.
\end{align*}
where $X_w$ is the usual composite of divided difference operators.
\end{ex}


\section{The twisted coproduct}\label{sec:coprod}
In this section we show that the structure algebra is a coalgebra with respect to some coproduct. This coproduct which we call a twisted coproduct is defined by iterating the reconstruction formula of the previous section.


\subsection{Coproduct on double Schubert classes}
Consider the product $S\otimes_{S^W} S \otimes_{S^W} S$ and denote by $g^{(i)}$ the respective element $g$ at the $i$th position, e.g., $g^{(2)}=1\otimes g \otimes 1$. Then the reconstruction formula for $u=v$ can be restated as follows
\[
\mu g^{(i)} = \sum_w \SC(u)_{I_w^{-1}}^{(ij)}\cdot   {}^u \tilde Y_{I_w}^{(j)}(g^{(j)}),\qquad i<j.
\]
In particular, we have
\[
\mu ({}^u \tilde Y_{I_w}^{(2)}(g^{(2)}))=\sum_x \SC(u)_{I_x^{-1}}^{(23)}\cdot   {}^u \tilde Y_{I_x}^{(3)} {}^u\tilde Y_{I_w}^{(3)}(g^{(3)}).
\]
Combining we get
\[
\mu^2g^{(1)}=\sum_{w} \SC(u)_{I_{w}^{-1}}^{(12)} \cdot \mu ({}^u \tilde Y_{I_w}^{(2)}(g^{(2)}))=\sum_{w,x} \SC(u)_{I_w^{-1}}^{(12)}\SC(u)_{I_x^{-1}}^{(23)} \cdot {}^u \tilde Y_{I_x}^{(3)} {}^u \tilde Y_{I_w}^{(3)}(g^{(3)}).
\]
On the other hand, we also have
\[
\mu g^{(1)}= \sum_{w} \SC(u)_{I_{w}^{-1}}^{(13)}\cdot  {}^u \tilde Y_{I_w}^{(3)}(g^{(3)}).
\]
Therefore, for any $g\in S$ we have
\begin{equation}\label{eq:mainrelY}
\mu \sum_{w} \SC(u)_{I_{w}^{-1}}^{(13)}\cdot  ({}^u \tilde Y_{I_w}(g))^{(3)}=
\sum_{x,y} \SC(u)_{I_x^{-1}}^{(12)}\SC(u)_{I_y^{-1}}^{(23)} \cdot ({}^u \tilde Y_{I_y} {}^u \tilde Y_{I_x} (g))^{(3)}.
\end{equation}

Recall that ${}^u \tilde Y_{I_w}=\sum_v d_{w,v}^u {}^u Y_{I_v}$ in $\mathbb{D}$. 
Since $\{ {}^u \tilde Y_{I_w}\}_{w}$ is a basis of $\mathbb{D}$, we can write
\begin{equation}\label{eq:Yrels}
{}^u \tilde Y_{I_y}  {}^u \tilde Y_{I_x}=\sum_z \nu(u)_{I_x^{-1},I_y^{-1}}^{I_z^{-1}} {}^u \tilde Y_{I_z},\qquad \text{ where }\nu(u)_{I_x,I_y}^{I_z} \in S.
\end{equation}

\begin{ex} 
In the additive case for $u=w_0$ we obtain ${}^{w_0}\tilde Y_{w}={}^{w_0}Y_{w_0w}=X_{ww_0}$. So
\[
X_{yw_0}X_{xw_0}=\sum_z \nu(w_0)_{x^{-1},y^{-1}}^{z^{-1}} X_{zw_0}
\]
which implies that 
\[
\nu(w_0)_{x,y}^z=\begin{cases}
1 & \text{ if }xw_0y=z\text{ and }l(x)+l(w_0y)=l(z)\\
0 & \text{ otherwise.}
\end{cases}
\]
\end{ex}

Combining \eqref{eq:mainrelY} and \eqref{eq:Yrels} we then obtain the following identity in $S\otimes_{S^W} S\otimes_{S^W}S$:
\begin{equation}  \label{eq:dscop}
\mu \SC(u)^{(13)}_{I_z} =\sum_{x,y} \SC(u)^{(12)}_{I_x}\SC(u)^{(23)}_{I_y} \nu(u)_{I_x,I_y}^{I_z}.
\end{equation}


\subsection{The twisted coproduct map}
Consider now the twisted tensor product $\SA\hat\otimes \SA$, where 
\[
(z_v)_v \hotimes (z_u')_u=(z_v\hotimes z_u')_{v,u}
\]
\[
v(s) z_v  \hotimes z_u'=z_v \hotimes sz'_u,\quad s\in S.
\]
In particular, we have
\[
(s\bullet \mathbf{1}) \hotimes \mathbf{1}= c(s) \hotimes \mathbf{1}= \mathbf{1}\hotimes s \mathbf{1}=\mathbf{1}\hotimes (s\odot \mathbf{1}).
\]
We may also look at this product as a product of bimodules, where $S$ acts on $\SA$ via the usual multiplication and via the characteristic map.

\begin{lem}
There is a commutative diagram
\[
\xymatrix{
S\otimes_{S^W} S \ar[r]^-{\varepsilon^{(14)}} \ar[d]_\rho & (S\otimes_{S^W} S)\otimes_S (S\otimes_{S^W} S) \ar[d]^{\rho\otimes \rho}\\
\SA \ar[r]^-{\De}& \SA\hotimes \SA,
}
\]
where $\varepsilon^{(14)}(a\otimes b)=a\otimes 1\otimes 1\otimes b$ and the map $\De\colon \SA \to  \SA\hotimes \SA$ is defined on twisted Schubert classes by 
\[
\De(\zeta(u)_{I_z})=\sum_{x,y} \zeta(u)_{I_x} \hotimes \zeta(u)_{I_y} c(\nu(u)_{I_x,I_y}^{I_z}).
\]

\end{lem}

\begin{proof} 
Consider the map $\phi\colon S\otimes_{S^W} S\otimes_{S^W} S \to (S\otimes_{S^W} S)\otimes_S (S\otimes_{S^W} S)$ defined by $a\otimes b\otimes c \mapsto (a\otimes b) \otimes (1 \otimes c)=(a\otimes 1) \otimes (b \otimes c)$.
By \eqref{eq:dscop} we get 
\begin{align*}\varepsilon^{(14)}\mu\SC(u)_{I_z} &=\phi(\mu\SC(u)_{I_z}^{(13)})\\
&=\sum_{x,y} \phi(\SC(u)^{(12)}_{I_x})\phi(\SC(u)^{(23)}_{I_y} {\nu(u)_{I_x,I_y}^{I_z}}^{(3)}))\\
&=\sum_{x,y}  \SC(u)_{I_x}  \otimes_S  (\SC(u)_{I_y} (1\otimes_{S^W} \nu(u)_{I_x,I_y}^{I_z})).\end{align*}
Applying $\rho\otimes \rho$ we obtain:
\[
(\rho\otimes\rho) (\varepsilon^{(14)}\SC(u)_{I_z})=\mu\sum_{x,y} \zeta(u)_{I_x}\hat \otimes \zeta(u)_{I_y}c(\nu(u)_{I_x,I_y}^{I_z}).
\]
The diagram then follows since $\SC(u)_{I_z}$ is a basis of $S\otimes_{S_W}S$ after inverting $\mu$.
\end{proof}

\begin{dfn} 
The map $\De$ of the lemma will be called a {\em twisted coproduct map}.
\end{dfn}

Observe that by the diagram, 
\[
\De(s_1 c(s_2))=\De(\rho(s_1\otimes s_2))=(\rho\otimes\rho)((s_1\otimes 1) \otimes (1\otimes s_2))=s_1\mathbf{1}\hotimes c(s_2).
\]
So $\De$ is a homomorphism of bimodules with respect to the usual $S$-action on the left factor and the $c$-twisted $S$-action on the right factor.
Moreover, since $\varepsilon^{(14)}$ is coassociative (i.e., we have $(\varepsilon^{(14)} \otimes \id) \circ \varepsilon^{(14)}=(\id \otimes \varepsilon^{(14)}) \circ \varepsilon^{(14)}$), so is $\De$.

Consider the multiplication map $S\otimes_{S^W} S \to S$ given by $s_1\otimes s_2\mapsto s_1s_2$. It defines the ring homomorphism $\varepsilon\colon \SA \to S$ by sending $\rho(s_1\otimes s_2)=s_1c(s_2)$ to $s_1s_2$. We call $\varepsilon$ the counit. The switch map $S\otimes_{S^W} S \to S\otimes_{S^W} S$ defines the map $\Upsilon\colon \SA \to \SA$ by $\rho(s_1\otimes s_2)\mapsto \rho(s_2\otimes s_1)$. Observe that by definition we have $\varepsilon\circ \Upsilon=\varepsilon$.

We then obtain
\[
(\varepsilon \hat \otimes \id)\circ \De (s_1c(s_2)) =(\varepsilon\hat \otimes \id)(\Upsilon(s_1\mathbf{1})\hotimes c(s_2))=(\varepsilon\hotimes \id) (c(s_1)\hotimes c(s_2))=s_1c(s_2)
\]
and, similarly,
\[
(\id \hat \otimes \varepsilon) \circ \De (s_1c(s_2)) =(\id \hat \otimes \varepsilon) (s_1\mathbf{1}\hotimes \Upsilon(c(s_2)))=(\id \hat \otimes \varepsilon)(s_1\hotimes s_2)=s_1c(s_2). \qedhere
\]
So there are counit diagrams:
\[
(\varepsilon\hat \otimes \id)\circ \De=(\id \hat \otimes \varepsilon)\circ \De=\id.
\]
Therefore, $(\De,\varepsilon)$ turns $\SA$ into a coalgebra. We call it the twisted coalgebra.


\section{Duoidal categories and structure algebras}\label{sec:duoidal}

The purpose of this section is to prove our main result that the structure algebra is a natural bimonoid in the duoidal category of $S$-bimodules (Theorem~\ref{thm:bimonoid}).
As a direct consequence of this fact we obtain that the double quotient of the structure algebra has a natural structure of a bialgebra (Corollary~\ref{cor:bialgeb})
and of a Hopf-algebroid (Corollary~\ref{cor:Hopfalg}).


\subsection{Duoidal categories and bimonoids}
Let $\mathcal{C}$ be a category with two monoidal structures $(\botimes,\bunit)$ and $(\hotimes,\hunit)$. 
Following \cite[\S6]{AM10} and \cite[\S1.2]{BCZ} $\mathcal{C}$ is called duoidal if there are
three morphisms 
\begin{equation}\label{eq:duodef2}
\iota\colon \hunit\botimes \hunit\longrightarrow \hunit,\qquad
\hat{\iota}\colon \bunit \longrightarrow \bunit \hotimes \bunit,\qquad 
\upsilon\colon \bunit\longrightarrow \hunit
\end{equation}
so that $(\hunit,\iota,\upsilon)$ is a monoid with respect to $\botimes$ and
$(\bunit,\hat\iota,\upsilon)$ is a comonoid with respect to $\hotimes$, and there is
a natural transform 
\begin{equation}\label{eq:duodef1}
\vartheta\colon (A\hotimes B)\botimes (C\hotimes D)\longrightarrow
(A\botimes C)\hotimes (B\botimes D)
\end{equation}
known as the \emph{interchange law} 
which satisfies the following diagrammatic-axioms:

\medskip 

\noindent
{\it Associativity}
\[
\compactstyle
\xymatrix{
\big((A_1\hotimes B_1)\botimes (A_2\hotimes B_2)\big)\botimes(A_3\hotimes B_3)\ar[r]^{\al}\ar[d]_{\vartheta\botimes\id}&
(A_1\hotimes B_1)\botimes \big((A_2\hotimes B_2)\botimes(A_3\hotimes B_3)\big)\ar[d]^{\id\botimes \vartheta}\\
\big((A_1\botimes A_2)\hotimes (B_1\botimes B_2)\big)\botimes(A_3\hotimes B_3)\ar[d]_{\vartheta}&
(A_1\hotimes B_1)\botimes
\big((A_2\botimes A_3)\hotimes (B_2\botimes B_3)\big)\ar[d]^{\vartheta}\\
\big((A_1\botimes A_2)\botimes A_3\big)\hotimes
\big((B_1\botimes B_2)\botimes B_3\big)\ar[r]^{\al\hotimes\al}&
\big(A_1\botimes (A_2\botimes A_3)\big)\hotimes
\big(B_1\botimes (B_2\botimes B_3)\big)
}
\]
\[
\compactstyle
\xymatrix{
\big((A_1\botimes B_1)\hotimes (A_2\botimes B_2)\big)\hotimes(A_3\botimes B_3)\ar[r]^{\al}\ar@{<-}[d]_{\vartheta\hotimes\id}&
(A_1\botimes B_1)\hotimes \big((A_2\botimes B_2)\hotimes(A_3\botimes B_3)\big)\ar@{<-}[d]^{\id\hotimes\vartheta}\\
\big((A_1\hotimes A_2)\botimes (B_1\hotimes B_2)\big)\hotimes(A_3\botimes B_3)\ar@{<-}[d]_{\vartheta}&
(A_1\botimes B_1)\hotimes
\big((A_2\hotimes A_3)\botimes (B_2\hotimes B_3)\big)\ar@{<-}[d]^{\vartheta}\\
\big((A_1\hotimes A_2)\hotimes A_3\big)\botimes
\big((B_1\hotimes B_2)\hotimes B_3\big)\ar[r]^{\al\botimes\al}&
\big(A_1\hotimes (A_2\hotimes A_3)\big)\botimes
\big(B_1\hotimes (B_2\hotimes B_3)\big)
},
\]
where $\al\colon (A \star B) \star C \to A\star (B\star C)$ and $\star=\botimes,\hotimes$.

\medskip

\noindent
{\it Unitality}
\[
\compactstyle
\xymatrix{
\bunit\botimes(A\hotimes B)\ar[r]^-{\hat\iota\botimes \id} & (\bunit\hotimes \bunit)\botimes(A\hotimes B) \ar[d]^{\vartheta}\\
A\hotimes B\ar[u]^{\wr} \ar[r]^-{\sim}& (\bunit\botimes A)\hotimes (\bunit\botimes B)}
\qquad
\xymatrix{
(A\hotimes B)\botimes \bunit\ar[r]^-{\id \botimes \hat\iota}&
(A\hotimes B)\botimes(\bunit\hotimes \bunit)\ar[d]^{\vartheta}\\
A\hotimes B\ar[u]^{\wr}\ar[r]^-{\sim}& (A\botimes \bunit)\hotimes (B\botimes \bunit)}
\]
\[
\compactstyle
\xymatrix{
\hunit\hotimes(A\botimes B)\ar@{<-}[r]^-{\iota\hotimes \id}&
(\hunit\botimes \hunit)\hotimes(A\botimes B)\ar@{<-}[d]^{\vartheta}\\
A\botimes B\ar@{<-}[u]^{\wr}\ar@{<-}[r]^-{\sim}& (\hunit\hotimes A)\botimes (\hunit\hotimes B)}
\qquad
\xymatrix{
(A\botimes B)\hotimes \hunit\ar@{<-}[r]^-{\id \hotimes \iota}&
(A\botimes B)\hotimes(\hunit\botimes \hunit)\ar@{<-}[d]^{\vartheta}\\
A\botimes B\ar@{<-}[u]^{\wr}\ar@{<-}[r]^-{\sim}& (A\hotimes \hunit)\botimes (B\hotimes \hunit)}.
\]
\noindent

A bimonoid over a duoidal category $\mathcal{C}$ is
an object $Z$ equipped with morphisms  
\[
\nabla\colon Z\botimes Z \longrightarrow Z,\qquad 
\rho\colon \bunit\longrightarrow Z,
\]
\[
\De\colon  Z \longrightarrow Z\hotimes Z,\qquad 
\varepsilon\colon Z\longrightarrow \hunit,
\]
such that $Z$ is a monoid with respect to
$(\nabla,\rho)$, 
$Z$ is a comonoid with respect to
$(\De,\varepsilon)$, and
 the following diagrams commute:
$$\xymatrix{
Z\botimes Z\ar[r]^{\nabla}\ar[d]_{\De\botimes\De} & Z\ar[r]^{\De} & Z\hotimes Z\\
(Z\hotimes Z)\botimes (Z\hotimes Z)
\ar[rr]^{\text{interchange}}_\vartheta&&
(Z\botimes Z)\hotimes (Z\botimes Z)
\ar[u]_{\nabla\hotimes\nabla}
}$$
$$
\xymatrix{
Z\botimes Z\ar[d]_{\nabla}\ar[r]^{\varepsilon\otimes\varepsilon}& \hunit\botimes \hunit\ar[d]^{\iota}\\
Z\ar[r]^{\varepsilon}& \hunit}
\qquad
\xymatrix{
Z\hotimes Z\ar@{<-}[d]_{\De}\ar@{<-}[r]^{\rho\hotimes\rho}& \bunit\hotimes \bunit\ar@{<-}[d]^{\hat\iota}\\
Z\ar@{<-}[r]^{\rho}& \bunit}
\qquad
\xymatrix@!=0.5pc{
&Z \ar[dr]^{\varepsilon}\\
\bunit\ar[ur]^{\rho}\ar[rr]^{\upsilon} && \hunit
}
$$


\subsection{Bimodules}
Let $R$ be a commutative ring and let $S$ be a commutative $R$-algebra. We denote by $\SBM$  the category of $S$-bimodules which are $R$-modules in the sense that $rx=xr$, $r\in R$, $x\in M$. 
Consider the following two tensor structures in $\SBM$:
\begin{itemize}
\item the usual bimodule tensor product
\begin{gather*}
M\botimes N : =M\otimes_R N\big/\big<s_1xs_2\otimes y=x\otimes s_1ys_2\mid s_1,s_2\in S,x\in M,y\in N\big>
\end{gather*}
\item the ordered (twisted) tensor product 
\begin{gather*}
M\hotimes N : =M\otimes_R N\big/\big<xs\otimes y=x\otimes sy\mid s\in S, x\in M,y\in N\big>
\end{gather*}
\end{itemize}
Observe that the unit $\hunit$ for $\hotimes$ is the bimodule $S$ itself, and the unit $\bunit$ for $\botimes$ is the bimodule $S\otimes_R S$. 
There is a natural  map 
\[
\vartheta\colon (A\hotimes B)\botimes (C\hotimes D)\longrightarrow
(A\botimes C)\hotimes (B\botimes D).
\]
We also have $
\iota \colon S\otimes S=S \stackrel{\id} \longrightarrow S$, the usual multiplication in $S$ gives
$\upsilon\colon S\otimes_R S \to S$, and there is a map
\[
\hat\iota\colon S\otimes_R S \longrightarrow \big(S\otimes_R S \big)\hotimes 
\big(S\otimes_R S\big) 
\]
given by $x{\otimes} y\mapsto (x{\otimes} 1){\otimes}(1{\otimes} y)$. 
According to \cite[Example 6.18]{AM10} (see also \cite[\S4.3]{BCZ}) this turns $\SBM$ into a duoidal category with
\[
(\botimes,\bunit) = (\botimes,S\otimes_R S)\quad\text{and}\quad
(\hotimes,\hunit) = (\hotimes,S). 
\]
Moreover, $Z=\bunit$ is a bimonoid in $\SBM$ with $\rho=\id$ and $\varepsilon=\upsilon$.


\subsection{Structure algebras} 
Consider now the $R$-algebra $S$, where $S$ is the FGR and the respective structure algebra $\SA$. 

On one side the twisted coproduct $\De\colon \SA\to \SA\hotimes \SA$ and the map $\varepsilon\colon \SA\to S=\hunit$ introduced before
turn $\SA$ into a comonoid in $\SBM$ with respect to $\hotimes$ ($\SA$ is the twisted coalgebra). On the other side, the usual multiplication on $\SA$ 
\[
\nabla\colon \SA\botimes \SA\longrightarrow \SA
\]
together with the Borel map $\rho\colon \bunit=S\otimes_{R}S \to \SA$ turn the structure algebra $\SA$ into a monoid in $\SBM$ with respect to $\botimes$, i.e., into the usual $S\otimes_{R} S$-algebra. 

Since $\bunit$ is a bimonoid and $\rho$ respects $\De$, we obtain our second main result:

\begin{thm}\label{thm:bimonoid}
The structure algebra $(\SA,\rho,\varepsilon)$ is a bimonoid in the duoidal category $\SBM$ with respect to the tensor structures $\botimes$ and $\hotimes$.
\end{thm}


\subsection{Base change and the double quotient}
Recall that a functor between two duoidal categories is called \emph{duoidal} if it is a (strong) monoid with respect to two monoidal structures and it commutes with \eqref{eq:duodef2} and \eqref{eq:duodef1} (see \cite[\S6]{AM10}). A natural transform between two duoidal functors is called \emph{duoidal} if it is monoidal with respect to two monoidal structures. The homomorphism $S\to T$ of $R$-algebras induces a duoidal functor by the base change
$\SBM \to T\text{-}\mathsf{BMod}_R$.
It follows by definition that this functor sends bimonoids to bimonoids. 

Assume there is an augmentation map $\ep\colon S\to R$. Let $I$ denote its kernel.
The functor 
$\SBM \to \RBM$
then sends a bimodule $M$ to the double quotient
$\widehat{M}=M/(IM+MI)$.
But in $\RBM$ two monoidal structures coincide, so a bimonoid turns into a bialgebra.
In particular, for the structure algebra we obtain:
\begin{cor}\label{cor:bialgeb}
The double quotient of the structure algebra $\widehat{\SA}=\SA/(I\SA +\SA I)$ is a bialgebra over $R$. 
\end{cor}


\subsection{The bi-antipode and the Hopf algebroid}
Since $\SA$ has a natural structure of a bimonoid in the duoidal category of bimodules, a natural question to ask is whether $\SA$ possesses a structure of a Hopf bimonoid. 
There are different ways of defining a Hopf bimonoid in duoidal categories. 
If we restrict ourselves to the duoidal category of bimodules $\SBM$, then the bimonoid $\SA$ is the Hopf bimonoid in the following sense (see \cite{BCZ}):
 
There is an $R$-linear map $\Upsilon \colon \SA\to \SA $ which we call a bi-antipode such that 
\begin{itemize}
\item $\Upsilon\circ \Upsilon = \id$,
\item $\Upsilon\circ \eta_l=\eta_r$ where $\eta_l\colon \hunit \to \SA$ (resp. $\eta_r\colon \hunit \to \SA$) is the $R$-linear  morphism defining the left (resp., right) $S$-module structure of the bimodule $\SA$,
\item there are commutative diagrams:
$$\compactstyle
\xymatrix{
\SA\ar[r]^{\varepsilon}\ar[d]_{\De}& \hunit \ar[r]^{\eta_l} & \SA\\
\SA\hotimes \SA\ar[rr]^{\id \otimes \Upsilon}&& \SA\otimes_r \SA\ar[u]_{\nabla}}\qquad
\xymatrix{
\SA\ar[r]^{\varepsilon}\ar[d]_{\De}& \hunit \ar[r]^{\eta_r} & \SA\\
\SA\hotimes \SA\ar[rr]^{\Upsilon\otimes \id}&& \SA\otimes_l \SA\ar[u]_{\nabla}}$$
\end{itemize}

$$\compactstyle
\xymatrix@!=0pc{
{\bunit}
    \ar[ddrrr]^{\rho}
    \ar[rrrr]^{\upsilon}
    \ar[ddd]_{\hat\iota}
&&&&
\hunit
    \ar@{=}[ddrrr]
    \ar[rrrr]^{\eta_l}
&&&&
{\bunit}
    \ar[ddrrr]^{\rho}
\\\\
&&&
\SA
    \ar[rrrr]^{\varepsilon}
    \ar[ddd]^{\De}
&&&&
\hunit
    \ar[rrrr]^{\eta_l}
&&&&
\SA\\
{\bunit}\hotimes{\bunit}
    \ar[ddrrr]^{\rho\hotimes \rho}
    \ar[rrrrrrrr]%
    |!{[urrr];[ddrrr]}\hole
&&&& &&&&
{\bunit}\otimes_r{\bunit}
    \ar[ddrrr]^{\rho\otimes_r \rho}
    \ar[uuu]^{\nabla}%
    |!{[ul];[urrr]}\hole
\\\\
&&& \SA\hotimes \SA
    \ar[rrrrrrrr]^{\id\otimes\Upsilon}
&&&& &&&&
\SA\otimes_r \SA
    \ar[uuu]_{\nabla}
}$$
and the similar diagram for $\eta_r$, $\otimes_l$ and $\Upsilon\otimes \id$.

Observe that if the two monoidal structures coincide ($\hotimes=\botimes$), then the above definition gives the diagrammatic definition of the usual antipode and, hence, of the Hopf algebra. Therefore, the previous corollary can be strengthen to:

\begin{cor}\label{cor:Hopfalg}
The double quotient of the structure algebra $\widehat{\SA}=\SA/(I\SA +\SA I)$ is a Hopf algebra over $R$. 
\end{cor}


\section{The dihedral case}

Following Examples~\ref{ex:dihedral} and \ref{ex:momentd} consider the dihedral root system of type $I_2(p)$, where $p$ is an odd prime. 
The purpose of this section is to compute the double quotient $\widehat\SA$ of the respective structure algebra in the additive case
or, informally speaking, to compute the virtual Hopf-ring $\CH(I_2(p);\Fp)$.

Recall that the coefficient ring $\mathcal{O}$ of $I_2(p)$ is the ring of integers
$\Z(\ga)$ in the formally real cyclotomic field $\Q(\ga)$,  where $\ga=\xi+\xi^{-1}$ and $\xi$ is the primitive $2p$-th root of unity. 
Consider the principal ideal 
$\mathfrak{p}=(\ga+2)=(\xi^{-1}(1+\xi)^2)$. For $0< k < p$ denote
\[
u_k:= \xi^k+\xi^{-k}+(-1)^k 2\quad \text{ and }\quad v_k:= \xi^k+\xi^{-k}-(-1)^k 2.
\]
Observe that $v_1=\ga+2$. We will use the following fact

\begin{lem}\label{Th:primeunitLemma}
The ideal $\mathfrak{p}$ is a prime ideal in $\mathcal{O}$ with residue field $\Fp$, $u_k$ is a unit, and $(v_k)=\mathfrak{p}$.
\end{lem}

\begin{proof} 
Since $\Q(\xi)/\Q(\ga)$ is a quadratic field extension, for any $x\in \Q(\ga)$ we get $N_\ga(x)^2=N_\xi(x)$ for the respective norms.
We have $N_\xi (\xi)=N_\xi(1+(-\xi)^k)=1$ and $N_\xi(1-(-\xi)^k)=p$. Therefore, we obtain
\[
N_\ga(u_k)^2 = N_\xi(\xi^{-k}(1+(-\xi)^k)^2) =1\quad \text{ and }\quad N_\ga(v_k)=p^2.
\]
Hence, $N_\ga(u_k)=\pm 1$ and $N_\ga(v_k)=\pm p$.
Since, $v_k/v_1\in \Z(\ga)$, the result follows.
\end{proof}

We may choose the positive and simple roots of $I_2(p)$ as
\[
\Phi_+=\{\al_1=1,\;\xi,\;\ldots,\; \xi^{p-2},\;\al_2=\xi^{p-1}\}
\]
(all roots have length 2), so that
the corresponding coroots are 
\[
\{\al_1^\vee=2,\; 2\xi,\; \ldots,\; \al_2^\vee=2\xi^{p-1}\}.
\]
The respective Coxeter group $W$ is generated by the simple reflections $s_1$ and $s_2$, 
and for each $0<k<p$ there are exactly two elements of length $k$ in $W$. Let
$\sigma_{k1}=\sigma_{\ldots s_2s_1}$ and $\sigma_{k2}=\sigma_{\ldots s_1s_2}$ denote the corresponding special Schubert classes.

By the Pieri formula~\cite[Ch.IV.(3.5)]{Hi82}, for each $0< k<p$ we have
\begin{align*}
c(\la)\cdot \sigma_{k-1,2}
& = \left<\la,\xi^0\right>\sigma_{k1}+
\left<\la,\xi^{p-k}\right>\sigma_{k2}\; \in \mathfrak{I}\\
c(\la)\cdot \sigma_{k-1,1} \in 
& = \left<\la,\xi^{k-1}\right>\sigma_{k1}+
\left<\la,\xi^{p-1}\right>\sigma_{k2}\; \in \mathfrak{I},
\end{align*}
where $\mathfrak{I}$ denotes the characteristic ideal in the augmented structure algebra $\overline{\SA}=\SA/I\SA$, that is the ideal generated by augmented elements in the image of the characteristic map.

Substituting $\la=2\xi^0$ and $\la=2\xi^{p-1}$, we get
\begin{align}
2\sigma_{k1}+
(\xi^{p-k}+\xi^{-(p-k)})\sigma_{k2}&\in \mathfrak{I},
\label{eq:I2Che1}\\
(\xi^{p-k}+\xi^{-(p-k)})\sigma_{k1}+
2\sigma_{k2}&\in \mathfrak{I}.
\label{eq:I2Che2}
\end{align}
Observe that the linear combinations of special Schubert classes \eqref{eq:I2Che1} and \eqref{eq:I2Che2} generate the degree $k$ component of $\mathfrak{I}$ as an abelian group. 

If $k<p$ is odd, by taking the sum of \eqref{eq:I2Che1} and \eqref{eq:I2Che2} we get 
\[
(\xi^{p-k}+\xi^{-(p-k)}+2)(\sigma_{k1}+\sigma_{k2})\in \mathfrak{I}.
\]
By Lemma~\ref{Th:primeunitLemma}, $\xi^{p-k}+\xi^{-(p-k)}+2$ is a unit in $\mathcal{O}$, hence, we have

\begin{equation}\label{eq:I2Che1new}
\sigma_{k1}+\sigma_{k2}\in \mathfrak{I}.
\end{equation}

Substituting it into \eqref{eq:I2Che2} we get 

\begin{equation}\label{eq:I2Che2new}
(\xi^{p-k}+\xi^{-(p-k)}-2)\sigma_{k1}\in \mathfrak{I}.
\end{equation}

Observe that elements \eqref{eq:I2Che1new} and \eqref{eq:I2Che2new} of the structure algebra 
generate the degree $k$ component of $\mathfrak{I}$ as an abelian group. 
This shows that for an odd $k$ we have
\[
\widehat\SA^{(k)}=\overline{\SA}/\mathfrak{I}^{(k)}\simeq \mathcal{O}/\mathfrak{p}=\Fp.
\]
with $pr(\sigma_{k1})=-pr(\sigma_{k2})$, where $pr$ denotes the quotient map $\SA\to \widehat\SA$.

Similarly, if $k<p$ is even, taking the difference of \eqref{eq:I2Che1} and \eqref{eq:I2Che2}, 
we obtain 
\[
\widehat\SA^{(k)}\simeq \mathcal{O}/\mathfrak{p}= \Fp,
\]
with 
$pr(\sigma_{k1})=pr(\sigma_{k2})$.

In the case $k=p$, by the Pieri formula we obtain
\begin{align*}
c(\la)\cdot \sigma_{p-1,2}
& = \left<\la,\xi^0\right>\sigma_{w_0}, \;\text{ and}\\
c(\la)\cdot \sigma_{p-1,1}
& = \left<\la,\xi^{p-1}\right>\sigma_{w_0}.
\end{align*} 
By similar arguments, substituting $\la=2\xi^0$ and $\la=2\xi^{p-1}$, we conclude that  $\widehat{Z}^{(p)}$ is trivial.
Summarizing the discussion, we obtain:

\begin{prop} 
The double quotient $\widehat\SA$ of the structure algebra $\SA$ of a dihedral root system $I_2(p)$, 
where $p$ is an odd prime, is isomorphic as a graded abelian group to
\[
\widehat\SA\simeq\; \stackrel{0}{\rule{0pc}{1pc}\mathcal{O}}\oplus \stackrel{1}{\rule{0pc}{1pc}\Fp}\oplus 
\stackrel{2}{\rule{0pc}{1pc}\Fp}\oplus\cdots\oplus \stackrel{p-1}{\rule{0pc}{1pc}\Fp}.
\]
\end{prop}

By one of our main results (Corollary~\ref{cor:Hopfalg}) the double quotient $\widehat{\SA}$ is the Hopf algebra.
By the Hopf--Borel theorem we have (over the finite field):
\[
\widehat\SA\otimes {\Fp} \simeq \Fp[x]/\left<x^p\right>,
\]
for some generator $x$.
We choose $x:=pr(\sigma_{s_1})=-pr(\sigma_{s_2})$. 
Observe that since $\widehat\SA$ is a graded Hopf algebra generated by one element, 
the coproduct and the antipode are given (over the finite field) by  
\[
\De(x) =x\otimes 1+1\otimes x,\qquad \Upsilon(x) = -x.
\]

On the other hand, the twisted coproduct can also be described using our formula: 
\[
\De(pr(\sigma_{w})) = 
\sum_{\begin{subarray}{c}
uv=w,\\
\ell(u)+\ell(v)=\ell(w)
\end{subarray}}
pr(\sigma_{u})\otimes pr(\sigma_{u}),\qquad
\Upsilon(pr(\sigma_{w}))  = (-1)^{\ell(w)} pr(\sigma_{w^{-1}}).
\]
Assume that $pr(\sigma_{k1}) = a_kx^k$.
Observe that 
\[
\De(x^k) = x^k\otimes 1+k x^{k-1}\otimes x+\cdots.
\]
On the other hand, 
\[
\De(\sigma_{k1}) = \sigma_{k1}\otimes 1+\sigma_{k-1,2}\otimes \sigma_{s_1}+\cdots
\]
This shows 
\[
a_kk \equiv (-1)^{k-1}a_{k-1}\mod p.
\]
As a result, 
\[
a_k \equiv \tfrac{(-1)^{k-1}}{k}a_{k-1}\equiv \cdots
\equiv \tfrac{(-1)^{k(k-1)/2}}{k!}\mod p.
\]
Combining all these together we obtain
\[
pr(\sigma_{k1})=\tfrac{(-1)^{k(k-1)/2}}{k!}x^k,
\]
and, since $pr(\sigma_{k1})=(-1)^k pr(\sigma_{k2})$, 
\[
pr(\sigma_{k2})=\tfrac{(-1)^{k(k+1)/2}}{k!}x^k.
\]


\bibliographystyle{alpha}

\end{document}